\documentclass[12pt,a4paper]{article}

\usepackage{amssymb}

\newcommand{\cont}{{\rm cont}}
\newcommand{\den}{{\rm den}}
\newcommand{\dens}{{\rm dens}}
\newcommand{\intl}{{\rm int}}

\newcommand{\loc}{{\rm loc}}
\newcommand{\vol}{{\rm vol}}
\newcommand{\mod}[1]{{\rm\ (mod\ }#1{\rm)}}
\newcommand{\CR}{\\\scriptstyle}
\newcommand{\cond}[1]
{{\renewcommand{\arraystretch}{.7}\renewcommand{\arraycolsep}{.2em}
\begin{array}{c}\scriptstyle#1\end{array}}}
\newcommand{\Cond}[1]
{{\makebox[0pt]{\renewcommand{\arraystretch}{.7}
$\begin{array}{c}\scriptstyle#1\end{array}$}}}

 \newcommand{\bm}[1]{\mbox{\boldmath $#1$}}
 \newcommand{\bms}[1]{\mbox{\scriptsize \boldmath $#1$}}

 \newcommand{\Adel}{\mathbb A}
 
 \newcommand{\NN}{\mathbb N}
 \newcommand{\RR}{\mathbb R}
 \newcommand{\CC}{\mathbb C}
 \newcommand{\ZZ}{\mathbb Z}

 \newcommand{\QQ}{\mathbb Q}

 \newcommand{\cS}{\mathcal{S}}

 \newcommand{\cF}{\mathcal{F}}

\newtheorem{theorem}{Theorem}
\newtheorem{lemma}{Lemma}
\newtheorem{prop}{Proposition}
\newtheorem{coro}{Corollary}

 \newcommand{\be}{\begin{equation}}
 \newcommand{\ee}{\end{equation}}

 \input{epsf}
 
\begin{document}
 \bibliographystyle{unsrt}

{ \parindent0pt

\begin{center}
{\large \bf Diffraction from visible lattice points \\
            and $k$th power free integers}

 \end{center}

 \vspace{5mm}
 
 \begin{center}  {\sc Michael Baake}$^1$, {\sc Robert V.~Moody}$^2$ 
                 and {\sc Peter A.~B.~Pleasants}$^3$
 
\vspace{5mm}

 1) Institut f\"ur Theoretische Physik, Universit\"at T\"ubingen, \\

 Auf der Morgenstelle 14, D-72076 T\"ubingen, Germany  \\

 2) Department of Mathematical Sciences, University of Alberta, \\

 Edmonton, Alberta T6G 2G1, Canada \\
 
 3) Department of Mathematics and Computing Science, \\
 University of the South Pacific, Suva, {}Fiji
 
\end{center}
 
\vspace{8mm}
\centerline{Dedicated to Ludwig Danzer on the occasion of his 70th birthday}
\vspace{25mm}
}

\begin{abstract} We prove that the set of visible points of any lattice
of dimension $n\geq 2$ has pure point diffraction spectrum,  
and we determine the diffraction spectrum explicitly.  This settles previous
speculation on the exact nature of the diffraction in this situation.
Using similar methods we show the same result for
the \mbox{1-dimensional} set of $k\/$th-power-free integers 
with $k\geq 2$. Of
special interest is the fact that neither of these sets is a Delone set
---  each has holes of unbounded inradius. We provide a careful
formulation of the mathematical ideas underlying the study of
diffraction from infinite point sets.

\end{abstract}

\clearpage

\subsection*{Introduction}

It has long been known that the diffraction spectrum of a crystal consists
of pure Bragg peaks only, being a pure point measure supported on the lattice 
dual to the lattice of periods.  
Until about 15 years ago it was tacitly believed that crystals
were the only discrete point sets with this property.  Now, however, we know
that many quasicrystals have pure point diffraction spectra (though in a
significantly different sense, since the locations of the peaks are no longer
discrete).  These quasicrystals are all Meyer sets, that is, sets $S\/$ that
are both uniformly discrete and relatively dense and whose difference sets 
$\Delta = S-S$ also have these properties. (For a discussion of Meyer sets,
see \cite{Moody}).  There has consequently been a feeling that a perfectly
diffracting discrete point set, if not precisely a Meyer set, must be 
closely related to a Meyer set. Indeed, the Meyer condition cannot possibly
be strictly necessary since, as we shall see later, adding
or removing a set of density zero does not alter the diffraction spectrum
of a discrete point set, and one can clearly destroy the relative denseness
of any uniformly discrete point set by removing a set of density zero.

In this paper we give some simple examples of perfectly diffractive
discrete point sets that deviate much further from the Meyer properties
than this; in fact we consider sets 
which, for arbitrarily large $D$, have a lattice of holes of inner diameter
at least $D$.  Such a set cannot differ from a Meyer set (or from a model
set, which is a particular case of a Meyer set) only by a set of density zero.
The sets comprising our examples are well known in number theory: they are the
sets of visible (or primitive) points of a lattice in any dimension $n\ge2$,
see \cite{Apostol,Hua} and the front cover of \cite{Apostol} for a picture
of the case $n=2$, and the $\mbox{1-dimensional}$ sets consisting of the
$k\/$th-power-free  integers for $k\ge2$, see \cite[\S 6.6]{Hua}.

As well as giving a rigorous derivation of the diffraction spectra that depends
on explicitly calculating the autocorrelation of the point set, we precede it
by a shorter derivation of the pure point parts of the spectra only, via the
{}Fourier transform of the point set. This takes a similar line to
previous attempts and uses a general result of A.~Hof. 
It provides a quick, elegant way of calculating the discrete parts of the
spectra which gives the correct results but does not have a full mathematical
justification at present. 

The rigorous approach is a response to the history of the problem for the 
visible lattice points, described in \cite{BGW}. 
In particular, there is a clear disagreement between earlier results
in \cite{Schroeder} and in \cite{Mosseri} regarding the nature
of the diffraction. This was partially resolved in favour
of \cite{Mosseri} by a calculation of the point part of the diffraction spectrum
in terms of Dirichlet series and its comparison with a real optical experiment
\cite{BGW}. The older numerical calculations in \cite{Schroeder}, using the fast
{}Fourier transform, suffer from an insufficient resolution and are misleading. 
In this article, we give the definitive description of the nature of the diffraction
spectrum by making the previous formal calculation of the pure point part
\cite{Mosseri,BGW} rigorous and by proving that there is no continuous part.

The paper is organized as follows. We start by describing the results
from number theory and related areas we need, then describe some
basic properties of the set of visible points.  Next we discuss the
background material of {}Fourier transforms and autocorrelations needed
for diffraction spectra.  This is first framed in the language of
tempered distributions and then connected to measure theoretical
considerations which are a basic part of the theory of diffraction. 
We then give the short intuitively motivated method of
computing the pure point part of the diffraction spectrum, both for the
visible points and the $k\/$th-power-free integers. 
The following two
sections are devoted to an explicit calculation of the autocorrelation
of the set of visible points, followed by a rigorous derivation of the
diffraction spectrum.  The ensuing section treats the set of
$k\/$th-power-free numbers in the same way, followed by an extensive outlook
and a summary.

Strictly speaking, neither the measure theoretical picture
nor the intuitive approach to computing the diffraction 
are necessary for the logic of the paper, and they could be omitted
by the reader whose primary concern is verifying the mathematical validity
of the results. 
But we would like to stress that the additional information provided
in this article is needed to locate our results in the wider context
of mathematical diffraction theory.

\subsection*{Tools from number theory}

Here we set out some results we need from number theory and related areas.

\subsubsection*{Notation}

We use the notation $(l,m)=\gcd(l,m)$ to denote the greatest common divisor of two
integers $l\/$ and $m$.  {}For integers $d\/$ and $x$, $d\mid x$ means that
$d\/$ is a divisor of $x$.  Summation conditions like ``$d\mid x$" are to
be interpreted as meaning that $d\/$ runs through positive divisors of
$x\/$ only (even when $x\/$ is negative).  The divisor function $\sigma(m)$
(defined for $m\in\ZZ^+$) counts the number of positive divisors of $m$,
so $\sigma(m) = \sum_{d|m} 1$.

Also, we will frequently use the $O$-notation for error estimates. {}For example,
we say that a function $f(r)$ is $O(1/r)$ if there is a constant $c$ such that 
$|f(r)|$ is bounded by $c/r$ for $r\geq 0$, see \cite[Sec.\ 3.2]{Apostol} for details.

\subsubsection*{Lattices}

A {\em lattice\/} in $\RR^n$ is a set $\Gamma$ of the form
\be
     \Gamma \; = \; \ZZ\bm{b}_1 \oplus \cdots \oplus \ZZ\bm{b}_n \, ,
\ee
where $\{\bm{b}_1, \ldots ,\bm{b}_n\}$ is a set of $n\/$ linearly
independent vectors called a {\em basis\/} of $\Gamma$.  We define
${\rm vol}(\Gamma)$ to be the volume of a fundamental region of the lattice,
e.g.\ of $\{t_1\bm b_1+\cdots+t_n\bm b_n\mid0\le t_1,\ldots,t_n<1\}$. 
Consequently, ${\rm vol}(\Gamma)$ can be calculated as 
${\rm vol}(\Gamma) = |\det( \bm{b}_1,\ldots,\bm{b}_n)|$
which turns out to be independent of the basis chosen.
Every lattice is uniformly discrete and relatively dense in $\RR^n$ and is a 
subgroup of $\RR^n$ under vector addition.

\begin{prop}\label{sphere} $\;$
   Let $\Gamma$ be a lattice and \bm a an arbitrary vector (not necessarily in
   $\Gamma$).  Let $N(R)$ be the number of points \bm x in $\Gamma+\bm a$ with
   $|\bm x|<R$.  Then there are constants $c_1$ and $c_2$, depending only
   on $\Gamma$, such that, for all $R>0$,
\be
  |N(R){\rm vol}(\Gamma) - v_nR^n| \; \le \; c_1R^{n-1}+c_2 \, ,
\ee
   where $v_n$ is the volume of the unit ball in $\RR^n$, 
   i.e.\ $v_n = \pi^{n/2} /\, \Gamma(1+{n\over 2})$.
\end{prop}
{\sc Proof}: The translates of the fundamental region of $\Gamma$ by vectors
in $\Gamma+\bm a$ tile $\RR^n$.  Let {\Large $v\/$} be the total volume of those
translates that meet the open ball $B_R(\bm{0})$ and $v\/$ the total volume of those
translates that lie entirely inside this ball.  Then the volume $v_nR^n$ of
the ball and the number $N(R){\rm vol}(\Gamma)$ are both bounded above by 
{\Large $v\/$} and below by $v$.  So their difference is at most 
{\Large $v\/$}$-v$, which is the total
volume of the translates of the fundamental region that meet the boundary of
the ball. If $D\/$ is the diameter\footnote{The diameter of a bounded set 
$S\subset\RR^n$ is the supremum of all distances between points of $S$.} 
of the fundamental region, then this volume is at most 
$v_n(R+D)^n-v_n(R-D)^n \le 2^{n}v_nDR^{n-1}$, when $R\ge D$, and
at most $v_n(R+D)^n<2^nv_nD^n$, when $R<D$.  The second of these 
estimates is obvious, while the first (when $R \ge D$) follows from
\begin{eqnarray}
(R+D)^n - (R-D)^n  & = & 
   \sum_{m=0}^{n} {n\choose m} R^{n-m} \left( D^m - (-D)^m \right) 
   \nonumber \\
& \le & 2 \sum_{m \; {\rm odd}} {n\choose m} R^{n-1} D \\
& = & 2 D R^{n-1} \sum_{m \; {\rm odd}} {n\choose m} 
\; = \; 2^n D R^{n-1} \, . \nonumber
\end{eqnarray}
This gives the result with
$c_1=2^{n}v_nD\/$ and $c_2=2^nv_nD^n$. \hfill $\square$

If $\Gamma$ is a lattice and $r\/$ is a nonzero real number then $r\Gamma$ is
also a lattice (a basis of $r\Gamma$ is $\{r\bm{b}_1, \ldots, r\bm{b}_n\}$).
When $r\/$ is an integer $r\Gamma\subseteq\Gamma$ and $r\Gamma$ is a
sublattice (and a subgroup) of $\Gamma$, of index $r^n$.  
{}For a nonzero integer $m\/$ and two
points $\bm a$ and $\bm b$ in $\Gamma$ we write $\bm a\equiv\bm b\mod{m\Gamma}$ to
mean that $\bm a-\bm b\in m\Gamma$.  With this notation we have the following
Chinese Remainder Theorem for a lattice $\Gamma$.
\begin{prop}\label{CRT} $\;$
  Let $\Gamma$ be a lattice and $\bm a^{}_1,\bm a^{}_2,\ldots,\bm a^{}_r\in\Gamma$.
  If $m^{}_1,m^{}_2,\ldots,m^{}_r\in\ZZ^+$ are chosen so that $(m^{}_j,m^{}_k)=1$, 
  for $1\le j<k\le r\/$, then there is a point $\bm a\in\Gamma$ such that
  the solutions of the simultaneous congruences
\be
  \bm x\equiv\bm a^{}_1\mod{m^{}_1\Gamma},\dots,\bm x\equiv\bm a^{}_r\mod{m^{}_r\Gamma}
\ee
  are precisely the points $\bm x\in\Gamma$ with
\be
   \bm x\equiv\bm a\mod{(m^{}_1\cdot m^{}_2\cdot\ldots\cdot m^{}_r)\Gamma} \, .
\ee
\end{prop}
{\sc Proof}: This follows by applying the Chinese Remainder Theorem
for integers \cite[Thm.~2.7.1]{Hua} to each coordinate with respect
to a basis of $\Gamma$. \hfill $\square$

We define the {\em content} of a nonzero lattice point $\bm x$ in a
lattice $\Gamma$ by
\be\label{cont}
     \cont(\bm x) \; := \; \max\{\, l\mid\bm x\in l\Gamma \} \, .
\ee
If $\bm x$ is expressed in terms of a basis of $\Gamma$, $\bm x = \sum x_j \bm{b}_j$,
then $\cont(\bm x) = \gcd(x_1,\ldots,x_n)$ (which is therefore
independent of the particular basis chosen).  {}For consistency and convenience,
we define $\cont(\bm 0)=\infty$.  {}For $m\in\ZZ^+$ and $\bm x\in\Gamma$, we have
\be
   \cont(m\bm x) \; = \; m\cdot\cont(\bm x) \, .
\ee
It is clear from (\ref{cont}) that for $\bm x\in\Gamma\setminus\{\bm0\}$
\be\label{contbound}
     \cont(\bm x) \; \le \; {|\bm x|\over L(\Gamma)} \, ,
\ee
where $|\bm x|$ is the Euclidean length of \bm x and $L(\Gamma)$ is the length
of the shortest nonzero vector in $\Gamma$.

The set $V=V(\Gamma)$ of {\em visible points} of a lattice $\Gamma$
(also known as the {\em primitive\/} points of $\Gamma$) is
\be                                    
     V \; := \; \{\bm x\in\Gamma\mid\cont(\bm x)=1 \} \, .
\ee
In terms of a lattice basis, $V\/$ consists of all points whose coordinates
have no common divisor.  These are precisely the lattice points that are
visible from the origin, in the sense that the line segment joining them to
the origin contains no other lattice point.  When $n=1$, $V\/$ consists of
two points, equidistant from 0, but otherwise $V\/$ is infinite and,
indeed, contains more than half the lattice points, as we shall see.
{}For $m\in\ZZ^+$ we have
\be                                    
     mV \; = \; \{\bm x\in\Gamma\mid\cont(\bm x)=m \} \, .
\ee

\subsubsection*{Density}

Let $S\/$ be a uniformly discrete set of points in $\RR^n$ (i.e.\ there is
a $c>0$ so that $\bm{x},\bm{y}\in S$ with $\bm{x}\neq \bm{y}$ implies 
$|\bm{x}-\bm{y}| \geq c$). We say $S\/$ has {\em natural density\/} 
${\cal D} = {\rm dens}(S)$ if
\be \label{density}
    {\cal D}_R \; = \; {\cal D}_R(S) \; := \; 
             \frac{|\{\bm x\in S\mid|\bm x|<R\}|}{v_n R^n} 
        \quad \longrightarrow \quad {\cal D}
\ee
as $R\to\infty$.   
The expression $|{\cal D} - {\cal D}_R|$ (a function of $R\/$) is called
the {\em error term\/} for the natural density.  {}For example, by
Proposition~\ref{sphere}, ${\cal D} = {\rm dens}(\Gamma)$ exists for every 
lattice $\Gamma$ and is equal to $1/{\rm vol}(\Gamma)$ (with error term $O(1/R)$).

If ${\rm dens}(S)$ exists and if $T\/$ is an orthogonal transformation then
${\rm dens}(T(S))={\rm dens}(S)$, but it is a failing of the natural density
that it is not always true that ${\rm dens}(T(S))={\rm dens}(S)/|\det(T)|$ when
$T\/$ is a general affine transformation, see Appendix.

When we use the word ``density" from now on we shall mean natural density.
 
\subsubsection*{Power-free numbers}

{}For an integer exponent $k\ge1$ the set $F=F_k$ of
{\em $k$th-power-free integers} is
\be                                    
     {}F \; := \; \{n\in\ZZ\mid
     \mbox{$n\/$ is not divisible by $d^k$ for any integer $d>1$} \} \, .
\ee
An equivalent characterization of the numbers $n\in {}F\/$ is that in their
prime power factorization $n=\pm p_1^{a_1}p_2^{a_2}\cdots p_r^{a_r}$ every
exponent $a_j$ is less than $k$.  In close analogy with $V$, the case $k=1$
is trivial, $F_1$ consisting  of just the two numbers $\pm1$, but for $k\ge 2$,
$F_k$ is infinite and contains more than half the integers. These numbers have
been studied for a long time, see \cite{Pillai,Mirsky} for two sources relevant
in our context.

\subsubsection*{Inclusion-exclusion}

The M\"obius function $\mu(m)$ is defined for $m\in\ZZ^+$ by
\be
\mu(m)  \; := \;
    \cases{1,     &when $m=1$,\cr
          (-1)^r, &when $m\/$ is a product of $r\/$ distinct primes,\cr
           0,     &when $m\/$ is divisible by the square of a prime.\cr}
\ee                     
It is multiplicative in the sense that $\mu(lm)=\mu(l)\mu(m)$ when
$(l,m)=1$ (and clearly $\mu(lm)=0$ when $(l,m)\ne1$).

A number of inversion formul\ae\ and variants of the inclusion-exclusion
principle can be expressed in terms of this function.  The result we need
in this paper is
\be
   \sum_{d\mid m}\mu(d) \; = \; \cases{1,&if $m=1$,\cr0,&if $m>1$,\cr}
\ee
see \cite[Thm.~6.3.1]{Hua}.  {}For a point \bm x in a lattice $\Gamma$
this shows that
\be\label{chiv}
   \chi^{}_V(\bm x) \; := \; \sum_{d\mid\cont(\bms x)}\mu(d)
\ee
is the characteristic function of the visible points $V\/$ of $\Gamma$
(except that it is undefined when $\bm x=\bm0$).

Similarly, since $d^k\mid x\/$ if and only if $d\/$ divides the
largest integer whose $k\/$th-power divides $x$,
\be\label{chik}
        \chi^{}_F(x) \; := \; \sum_{d^k\mid x}\mu(d)
\ee
is the characteristic function of the $k\/$th-power-free integers
(except for being undefined when $x=0$).

\subsubsection*{Dirichlet series}

The well-known Riemann zeta-function (see, for example, \cite[Ch.\ 12]{Apostol} 
or \cite[\S9.2]{Hua}) is defined for complex numbers $s\/$ with ${\rm Re}(s)>1$ by
\be\label{zeta}
       \zeta(s) \; := \; \sum_{m=1}^\infty{1\over m^s}
                \; = \; \prod_p\left(1-{1\over p^s}\right)^{-1},
\ee
where the sum is called a {\em Dirichlet series\/} and the product (in which
$p\/$ runs through all positive prime numbers) is called an {\em Euler
product}.  The sum and product are absolutely convergent for $s\/$ in the
half-plane ${\rm Re}(s) >1$, but $\zeta(s)$ can be meromorphically continued to
the whole of $\CC$, for example $\zeta(0)=-1/2$.  The only singularity of
$\zeta(s)$ is a simple pole at $s=1$ with residue 1.  Using the Euler product,
we see that the Dirichlet series of $1/\zeta(s)$ is given by
\be\label{invzeta}
       {1\over\zeta(s)} \; = \; \prod_p\left(1-{1\over p^s}\right)
                        \; = \; \sum_{m=1}^\infty{\mu(m)\over m^s} \, .
\ee
This function has infinitely many poles (all in the half-plane
${\rm Re}(s) <1$).  Its value at 1 is 0 and its value at 0 is $-2$.

Another Dirichlet series we shall encounter is
\be\label{xi}
    \xi(s) \; := \; \sum_{m=1}^\infty{\mu(m)\sigma(m)\over m^s}
           \; = \;  \prod_p\left(1-{2\over p^s}\right) .
\ee
Again, both the Dirichlet series and Euler product are absolutely
convergent in the half-plane ${\rm Re}(s) >1$.  It can also be seen from the
Euler product that $\xi(1)=0$.

\subsection*{Visible points of a lattice}

Here we summarize some elementary and well-known properties of the set
of visible points $V\/$ of a lattice $\Gamma$, together with complete proofs.

Since $\Gamma$ is a free Abelian group of rank $n$, its automorphism
group, $\mbox{Aut}(\Gamma)$, is isomorphic to the matrix group $GL(n,\ZZ)$.
Explicit isomorphisms can be found by taking coordinates with respect
to any lattice basis.
\begin{prop} \label{symm} $\;$
  The orbits of the action of $GL(n,\ZZ)$ on $\Gamma$ are the sets
  $m V$, $m\in\NN_0$. In particular, $GL(n,\ZZ)$ acts transitively on $V$.
\end{prop}
{\sc Proof}: Since the elements of $GL(n,\ZZ)$ cannot decrease content and
are invertible, they preserve content.  Hence each set $mV\/$ is invariant
under $GL(n,\ZZ)$.  The transitivity of $GL(n,\ZZ)$ on $V\/$ can be seen
from the facts that every visible point belongs to some basis of $\Gamma$
\cite[\S3,Thm.~5]{GL} and that any two bases of $\Gamma$ are related by a
transformation in $GL(n,\ZZ)$.  {}For $m\in\ZZ^+$ the transitivity of
$GL(n,\ZZ)$ on $mV\/$ follows from its transitivity on $V$.  Transitivity
on the singleton orbit $0\cdot V=\{\bm0\}$ is trivial. \hfill $\square$

\begin{prop} $\;$
   $V\/$ is uniformly discrete, but has arbitrarily large holes.  Moreover,
   for any $r>0$, there is a set of holes in $V\/$ of inradius at least $r\/$
   whose centres have positive density.
\end{prop}
{\sc Proof}: The uniform discreteness is trivial, as $V\/$ is a subset of a
lattice.  Now let $C=\{\bm a^{}_1,\dots,\bm a_s\}$ be any finite configuration
of points in $\Gamma$ (e.g.\ all points in a ball or a cube).  
Choose $s\/$ integer moduli $m_1,\ldots,m_s$ that are
$>1$ and coprime in pairs (for example, they could be the first $s\/$ primes). 
By Proposition~\ref{CRT} there is a point $\bm a\in \Gamma$ with
\be
  \bm a\equiv-\bm a^{}_1\mod{m^{}_1\Gamma} \, ,
       \dots, \, \bm a\equiv-\bm a^{}_s\mod{m^{}_s\Gamma} \, .
\ee
Now for any $\bm x\equiv\bm a\mod{m_1m_2\ldots m_s\Gamma}$ the configuration
$C+\bm x=\{\bm a^{}_1+\bm x,\dots,\bm a_s+\bm x\}$ is congruent, in the geometric
sense, to $C\/$ but no point in $C+\bm x$ is visible, since $\bm a_j+\bm x\in
m_j\Gamma$.  The points \bm x have density 
$\dens(\Gamma)/(m^{}_1m^{}_2\ldots m_s)>0$. \hfill $\square$

{}For $n=2$, and with the density of holes not mentioned, this is Thm.~5.29 of
\cite{Apostol}.  We note that the hole nearest the origin
provided by this argument can be expected to be at a distance of the order
$s^s$ and that the density guaranteed for holes of inradius $r\/$ is of the
order $r^{-n^2 r^n}$, so large holes, while having positive density, are
probably extremely sparse.

This proposition shows that $V$, though uniformly discrete, is not relatively
dense, and hence not a Delone set. Consequently it is not a Meyer set either.
(Recall that $\Lambda$ is a Meyer set if and only if both $\Lambda$ and
$\Lambda-\Lambda := \{ x-y\mid x,y\in\Lambda\}$ are Delone sets \cite{Lagarias}.)  
Also, $V$ cannot be transformed into a Delone set by adding a set of zero density. 
However, we do have:
\begin{prop} \label{diffprop} $\;$
     If $n\ge2$ then $V-V \; = \; \Gamma$.
\end{prop}
{\sc Proof}: Clearly, $V-V\subseteq\Gamma$.
Now let $\bm x = \sum_{j=1}^{n} x_j \bm{b}_j\in\Gamma$, where $n\geq 2$.
Then, $(x_1,x_2,\ldots,x_n)=(x_1+1,1,x_3,\ldots,x_n)-(1,1-x_2,0,\ldots,0)$,
and \bm x is the difference of two visible points. \hfill $\square$

The existence of arbitrarily large holes also implies that the set of
visible points cannot have a uniform density (not, at least, when its
density is positive in some sense, as is the case for $n\ge2$).
It does have a natural density, however:
\begin{prop} \label{freq1} $\;$
  The visible points $V\/$ of a lattice $\Gamma\in\RR^n$ have a natural
  density given by
\be
        {\rm dens}(V) \; = \; \frac{{\rm dens}(\Gamma)}{\zeta(n)} \, ,
\ee
with error term $O(1/R)$ when $n\ne2$ and $O(\log R/R)$ when $n=2$.
\end{prop}
This is a standard example of the use of M\"obius inversion
given (at least for the case $\Gamma=\ZZ^2$) in most introductory
number theory books (for example, \cite[Thm.~6.6.3]{Hua} and
\cite[Thm.~3.9]{Apostol}).  In these particular references the averages
are taken over triangles and squares, respectively, instead of balls.
Indeed the density is independent of the shape of the region 
averaged over. We say more about this in the Appendix.

{\sc Proof}:
The proposition is trivially true for $n=1$, when the pole of $\zeta(s)$
at 1 gives a density of 0, so we assume from now on that $n\ge2$.

The density of $V\/$, if it exists, is the limit as $R\to\infty$ of
\be\label{Dvis}
  {1\over v_nR^n}\sum_\Cond{|\bms x|< R\CR\bms x\in V}1 \, ,
\ee
which by (\ref{chiv}) is
\be\label{sumDvis}
  {1\over v_nR^n}\sum_\cond{\bms x\in\Gamma\setminus\{\bms0\}\CR|\bms x|< R}\quad
  \sum_\Cond{m\mid\cont(\bms x)}\mu(m)
  \; = \; {1\over v_nR^n}\sum_\cond{1\le m<cR}\mu(m) \;
  \sum_\Cond{\bms x\in m\Gamma\setminus\{\bms 0\}\CR|\bms x|< R}\;1 \; ,
\ee
where $c=1/L(\Gamma)$ (using (\ref{contbound})).  The inner sum is equal
to the number of nonzero points $\bm y\in\Gamma$ with $|\bm y|<R/m$,
which by Proposition~\ref{sphere} is
\be
  {v_n\over{\rm vol}(\Gamma)}\left({R\over m}\right)^n
  +O\left(\left({R\over m}\right)^{n-1}\right)+O(1) \, .
\ee
Substituting this into the right hand side of (\ref{sumDvis}) gives
\be
{\rm dens}(\Gamma)\!\!\sum_\cond{1\le m<cR}\!{\mu(m)\over m^n}
+O\biggl({1\over R}\sum_\cond{1\le m<cR}{1\over m^{n-1}}\biggr)
+O\left({1\over R^{n-1}}\right)
\ee
which tends to ${\rm dens}(\Gamma)/\zeta(n)$ as $R\to\infty$ when $n\ge2$ by
(\ref{invzeta}).  The total error is $O(1/R^{n-1})$, from the last term and
the tail of the sum in the main term, and $O(1/R)$ (or $O(\log R/R)$ when $n=2$
and the series diverges logarithmically), from the middle term. \hfill $\square$

Calculations of densities by M\"obius inversion form the core of this paper.
This is the first of many.

\subsection*{Diffraction spectra}

In this section we assemble the facts we need about distributions, {}Fourier
transforms and diffraction spectra. The mathematics underlying diffraction
is quite subtle and needs to be spelt out carefully. Although the 
discussion in this section does not contain
much that is new, it is nonetheless difficult to extract it all
from any convenient source. We shall use a formulation based upon
tempered distributions. {}For an essentially parallel approach which starts 
with measures, we refer to \cite{Hof}. We shall link the two approaches
later in this section.

\subsubsection*{Autocorrelations}

In dealing with diffraction, and therefore {}Fourier transforms, it is
appropriate to use tempered distributions, whose test space $\cal S$
consists of the Schwartz functions (also known as ``rapidly decreasing
functions"). We refer to (\cite{Schwartz,Rudin}) for this and the details 
on the standard topology which is used to describe convergence in $\cal S$.  
 A {\em tempered distribution} $T\/$ is a continuous (in the
sense of this topology) linear functional on $\cal S$. The space of
tempered distributions is denoted by ${\cal S}'$ and is equipped with the 
weak*-topology.  A simple example is
$\delta_{\bms x}$, {\em Dirac's delta-distribution}\footnote{Nowadays 
it is usually called Dirac's point measure for reasons that will become 
clear shortly.} at the point $\bm x$,
defined by $(\delta_{\bms x},\psi) := \psi(\bm x)$ for all $\psi\in\cal S$.
We prefer this notation to $\delta_{\bms x}(\psi)$ because it emphasizes
the duality between ${\cal S}'$ and ${\cal S}$.

If $S\/$ is a uniformly discrete subset of $\RR^n$, we call
\be\label{comb}
     \omega^{}_S \; := \; \sum_{\bms{x}\in S} \delta_{\bms x}
\ee
its {\em Dirac comb}.  It is also a tempered distribution, and in fact the
sum (\ref{comb}) is convergent (with any ordering) in the weak*-topology.

To describe the diffraction from a Dirac comb $\omega = \omega^{}_S$ we
need its {\em natural autocorrelation distribution} (also called its
generalized Patterson function) defined by
\be\label{gamma}
     \gamma^{}_{\omega} \; := \; \lim_{R\to\infty} {1\over v_nR^n}
        \sum_{\bms x,\bms y \in S_R} \delta_{\bms x - \bms y}\,,
\ee
where $S_R = S \cap B_R(\bm 0)$ and $B_R(\bm 0)$ is the ball of radius $R\/$ 
centre $\bm 0$. We shall simply call $\gamma^{}_{\omega}$ the ``autocorrelation" 
of $S\/$ from now on.  The existence of this limit in the weak*-topology is
a prerequisite for the diffraction spectrum to be well defined. (The word
``natural" refers to the use of the expanding ball $B_R(\bm{0})$ 
in the averaging process. Replacing it
by an expanding region of some other shape might lead to a different limit.)
It is clear from the definition that enlarging or diminishing $S\/$ by a
set of density~0 does not change its autocorrelation $\gamma^{}_{\omega}$.
In particular, adding or removing any finite number of points does not
change $\gamma^{}_{\omega}$, as for natural density.

We say that $S\/$ has {\em finite local complexity} if 
$\Delta = S-S$  is discrete and closed in $\RR^n$
(as is clearly the case for all sets whose diffraction
spectra we seek in this paper). Such an $S\/$ is 
uniformly discrete because $\bm{0}\in\Delta$ is isolated. 
The existence of the limit (\ref{gamma}) in the
weak*-topology is now controlled by the following result.
\begin{lemma}\label{ACexist} $\;$
 Let $S\/$ be a set of finite local complexity and $\omega=\omega_S$ its
 Dirac comb.  Then $S\/$ has a natural autocorrelation $\gamma^{}_{\omega}$
 if and only if the coefficients $w({\bm a})$ 
\be\label{coeff}
   w(\bm a) \; := \; \lim_{R\to\infty}{1\over v_nR^n}\!\!
\sum_\cond{|\bms x|,|\bms x-\bms a|< R\CR\bms x,\bms x-\bms a\in S}
\!\!\!\!\!\!1
\ee
 exist, i.e.\ the right hand side is convergent for all $\bm a\in\Delta$.  
 In this case, $w(\bm a) \ge 0$ for all $\bm a \in\Delta$ and 
 $\gamma^{}_{\omega}$ is the tempered distribution of 
 positive type\footnote{See next section for a definition.} given by 
 the weak*-convergent sum
\be\label{Auto}
  \gamma^{}_{\omega} \; = \;
            \sum_{\bms a \in \Delta} w(\bm a)\delta_{\bms a} \, .
\ee
\end{lemma}
{\sc Proof}: The existence of $\gamma^{}_{\omega}$ clearly implies the
existence of $w(\bm a)$ for all $\bm a \in \Delta$ because $\Delta$ is
discrete by assumption and the Schwartz space $\cal S$ contains all 
$C^{\infty}$-functions of compact support, whence we can focus on any 
individual $\bm a$.

Conversely, assume that the $w(\bm a)$ exist (they are then clearly $\ge 0$).
Since $\Delta$ is closed and discrete, its intersection with any compact
subset of $\RR^n$ contains only finitely many points. Consequently,
the $w(\bm a)$ are locally summable and the right hand side of (\ref{Auto})
defines a distribution over the space $\cal D$ of all $C^{\infty}$-functions of
compact support. We have to show that it is actually also a tempered distribution.
This follows from the translation boundedness\footnote{We will explain this
below in the context of measures.} of $\omega$ which is then inherited by
$\gamma^{}_{\omega}$, see \cite[Prop.\ 2.2]{Hof}. {}Finally, $\gamma^{}_{\omega}$
can also be written as a certain volume-normalized convolution (see below)
which implies that it is a distribution of positive type. \hfill $\square$

\subsubsection*{Fourier transforms}

The {}Fourier transform $\hat{T}$ of a tempered distribution $T\in\cal S'$
is defined by $(\hat{T},\psi) := (T,\hat{\psi})$, where we use the definition
\be
   \hat{\psi}(\bm y) \; := \; 
          \int_{\RR^n} e^{-2\pi i \bms y\cdot \bms x}\, \psi(\bm x) d \bm x 
\ee 
for the {}Fourier transform of functions $\psi\in{\cal S}$.  The {}Fourier
transform maps the space $\cal S$ onto itself and is continuous on $\cal S$
in the standard topology for Schwartz functions \cite[Thm.~7.7]{Rudin}, hence
it maps ${\cal S}'$ onto itself and is continuous on ${\cal S}'$ in the
weak*-topology \cite[Thm.~7.15]{Rudin}.

As special cases, we mention $\hat\delta_{\bms 0}=1$ and the well-known
{\em Poisson summation formula} for lattice Dirac combs
\be \label{poisson}
\hat{\omega}^{}_{\Gamma} \; = \; {\rm dens}(\Gamma) \cdot \omega^{}_{\Gamma^*}
\ee
where $\Gamma^*$ is the dual or reciprocal lattice defined by
\be
\Gamma^* \; := \;
   \{\bm y\mid\bm y\cdot\bm x \in\ZZ\,,\;\mbox{ for all }\bm x\in\Gamma\}\,.
\ee
(Eq.~(\ref{poisson}) can easily be derived from Poisson's summation formula 
for Schwartz functions \cite[p.\ 254]{Schwartz}.)

{}For a set $S$ of finite local complexity with autocorrelation $\gamma^{}_{\omega}$, 
its  {\em diffraction pattern\/} (also called its {\em diffraction distribution\/} or 
{\em diffraction spectrum\/}) is the {}Fourier transform $\hat{\gamma}^{}_{\omega}$. 
In view of the remarks above, $\hat{\gamma}^{}_{\omega}$ 
is a tempered distribution. It is also a positive measure, as we shall see.

The autocorrelation of a lattice $\Gamma$, for example,
is supported on $\Gamma$ itself (since $\Delta = \Gamma - \Gamma = \Gamma\/$)
and each peak has equal amplitude ${\rm dens}(\Gamma)$.  So the autocorrelation
is $\gamma^{}_{\omega}={\rm dens}(\Gamma)\omega^{}_\Gamma$ and the
corresponding diffraction spectrum is 
$\hat{\gamma}^{}_{\omega} = {\rm dens}(\Gamma)^2\omega^{}_{\Gamma^*}$, 
a constant multiple of the {}Fourier transform of $\omega^{}_\Gamma$ itself.
However, as we shall see, the {}Fourier transform of a general point set (even
of finite local complexity) does not describe
its diffraction in such a simple way.

\subsubsection*{Pure point distributions}

Our principal concern is with showing that the visible points and
$k${\em th}-power-free points have a pure point spectrum. In this section we
consider a special class of point measures which we will use in the 
sequel. We already borrow from the terminology of measures here, although
we will establish the precise connection only in the next section.

Consider an arbitrary complex point measure. It can be expressed in the form
\be\label{pointDist}
     \nu \; = \; \sum_{\bms x\in S}w(\bm x)\delta_{\bms x} \, ,
\ee
where the point set $S\/$ is countable, but not necessarily uniformly discrete,
and the coefficients or {\em weights} $w(\bm x) \in \CC$ are not necessarily 
constant. Note that the weights may be complex numbers.
Let us assume in addition that the measure is translation bounded.
This can be expressed as the condition that for
every compact set $K\subset\RR^n$ the sum
\be \label{transBoundedPointMeasure}
   \sum_{\bms x\in S \cap (K+\bms a)}|w(\bm x)|
\ee
is convergent and bounded uniformly in \bm a.  
 We denote the space of all translation bounded point measures
by ${\cal T}$. All these measures are tempered, and we identify
$\nu$ with $T_{\nu}$, the corresponding tempered distribution.

As a subset of ${\cal S}'$, $\cal T$ is not closed in the weak*-topology.
{}For example, the sequence of pure point distributions
\be\label{integral}
       \{j^{-n}\omega^{}_{\ZZ^n/j}\}
\ee
tends to the constant function 1 as $j\to\infty$.  (For a test function $\psi$
the numbers $j^{-n}\omega^{}_{\ZZ^n/j}(\psi)$ are approximating sums to the
integral\footnote{Seen as a sequence of measures, (\ref{integral})
weak*-converges to Lebesgue measure.} of $\psi$.)  
A similar argument shows that every bounded continuous
function is a limit of pure point distributions, compare the more
detailed discussion in \cite[Sec.\ IV.5]{RS}. Unfortunately, as this example shows,
the weak*-limit of such a sequence is not, in general, the same as its
pointwise limit.

However, taking pointwise limits in ${\cal T}$ is justifiable under certain 
circumstances, and fortunately these apply in the cases of interest to us here.
We introduce a ``locally defined'' norm on $\cal T$ by
\be
   \|\nu\|^{}_\loc \; := \; \sup_K \int_K d|\nu| \; = \;
   \sup_K\sum_{\bms x\in S\cap K}|w(\bm x)| \, ,
\ee
where the supremum is taken over all compact sets $K\/$ of diameter $< 1$.
This norm defines a topology on $\cal T$ stronger than the weak*-topology and
it provides a simultaneous ``$M\/$-test" for pointwise and weak*-convergence
of infinite sums of translation bounded point measures:
\begin{lemma}\label{Mtest} $\;$
  If $\nu_j\in\mathcal{T}$ for $j\in\ZZ^+$ and $\sum_{j=1}^\infty\|\nu_j\|_\loc$
  is convergent then $\sum_{j=1}^\infty\nu_j$ is pointwise and
  weak*-convergent to the same sum $\nu\in\mathcal{T}$.
\end{lemma}

{\sc Proof}: Let $\nu_j=\sum_{\bms x\in S_j}w_j(\bm x)\delta_{\bms x}$
and choose a fixed covering of $\{K_m\}_{m\in\ZZ^+}$ of $\RR^n$ by
compact sets of diameter $<1$.  For any $\psi\in\cS$ we have
\begin{eqnarray}  \label{nujpsi}
   |(\nu_j,\psi)| & = &
   \Bigl|\sum_{\bms x\in S_j}w_j(\bm x)\psi(\bm x)\Bigr|
   \; \le \; 
   \sum_{m=1}^\infty \; \sum_{\bms x\in S_j\cap K_m}|w_j(\bm x)\psi(\bm x)|
   \nonumber \\
   & \le &
   \|\nu_j\|_\loc\sum_{m=1}^\infty\|\psi^{}_{K_m}\| 
   \; = \; C_{\psi} \cdot \|\nu_j\|_\loc \, ,
\end{eqnarray}
where $\psi^{}_{K_m}$ is the restriction of $\psi$ to $K_m$ and $C_{\psi}<\infty$ 
is a constant that depends only on $\psi$ and the covering chosen.
Hence $\sum_{j=1}^\infty(\nu_j,\psi)$ is absolutely convergent
(by comparison with $\sum\|\nu_j\|_\loc$) and so $\sum\nu_j$ is
weak*-convergent to a distribution $\nu$.

Also
\be
   (\nu,\psi) \; = \;
   \sum_{j=1}^\infty(\nu_j,\psi) \; = \;
   \sum_{j=1}^\infty\sum_{\bms x\in S_j}w_j(\bm x)\psi(\bm x) \; = \;
   \sum_{\bms x\in\bigcup S_j}\sum_{j=1}^\infty w_j(\bm x)\psi(\bm x) \, ,
\ee
the reversal of the order of summation (with $w_j(\bm x):=0$ whenever
$\bm x \not\in S_j$) being justified by the fact that
the double sum is absolutely convergent, in view of (\ref{nujpsi}).  Hence
\be
   \nu \; = \; 
   \sum_{\bms x\in\bigcup S_j}\sum_{j=1}^\infty w_j(\bm x)\delta_{\bms x}
\ee
is the pointwise sum of the $\nu_j$'s.

Finally, if $K\subseteq\RR^n$ is any compact set with diameter $<1$ and we
write $S:=\bigcup S_j$ and $w(\bm x):=\sum_{j=1}^\infty w_j(\bm x)$ then
\be
   \sum_{\bms x\in S\cap(K+\bms a)}|w(\bm x)| \; \le \;
   \sum_{j=1}^\infty\sum_{\bms x\in S_j\cap(K+\bms a)}|w_j(\bm x)| \; \le \;
   \sum_{j=1}^\infty\|\nu_j\|_\loc \, .
\ee
Hence $\nu$ is translation bounded, and therefore tempered.    \hfill $\square$

Note that if $\mu_j$ is the distribution in (\ref{integral}) then
$\mu_{j+1}-\mu_j$ is translation bounded but
\be
    \|\mu_{j+1}-\mu_j\|_\loc  \; \ge  \; 
     {v_n\over2^{n-1}}-{1\over j(j+1)} \, ,
\ee
where $v_n$ is the volume of the $n\/$-dimensional unit ball.  (The supports
of $\mu_j$ and $\mu_{j+1}$ intersect in the integer points only and there
is almost no cancellation in calculating the norm.)  So the norms of the
differences are bounded below by a positive constant for large $j\/$ and
their sum diverges to infinity, thus failing to satisfy the hypothesis of
the lemma.  (An easier calculation along the same lines is
$\|\mu_{2^{k+1}}-\mu_{2^{k}}\|\ge v_n/2^n$.)

\subsubsection*{Distributions and measures}

The diffraction pattern $\hat{\gamma}^{}_{\omega}$
is a tempered distribution. However, it is also a positive measure. 
This remarkable fact is indispensible for the general theory of 
diffraction, making available to it a vast array of concepts and
tools. {}For example, it makes it immediately evident that the diffraction
pattern may be viewed as having a pure point part and a continuous
part, a fact of considerable physical significance.

The next two subsections sketch out this distribution theory --
measure theory connection in the context of diffraction.
The proofs of the main theorems of this paper in no way depend on
this background material, which can therefore be ignored as far as
verifying the results is concerned.  Nevertheless, it is desirable
to realize the natural connection to measure theory, both to see our
results in their proper setting and because it is this picture
that has a generalization to the diffraction theory of translation bounded 
measures on locally compact Abelian groups which is an appropriate setting 
for more general questions, compare \cite{HR,Arga}. 

A distribution $T$ is {\em positive} if for all positive test
functions $\psi$, $(T,\psi) \ge 0$. A distribution $T$ is of 
{\em positive type} or is  {\em positive definite} 
if for all test functions $\psi$, \, $(T,\psi * \psi^*) \ge 0$,
where $\psi^*(x) := \overline{\psi(-x)}$.
These two concepts are related by the Bochner-Schwarz
theorem (see \cite[Thm.\ VII.XVIII]{Schwartz} or \cite[Thm.\ IX.10]{RS})
which asserts that a tempered distribution
is positive if and only if it is the {}Fourier transform of a tempered
distribution of positive type.

Let ${\cal C}$ denote the space of complex-valued continuous functions
of compact support on $\RR^n$ and let $\| \cdot \|$ denote the 
supremum norm on ${\cal C}$. A (complex) {\em measure} $\nu$ on $\RR^n$ is 
defined as a linear functional on ${\cal C}$ such that
for every compact set $K\subset\RR^n $ there is a constant $a^{}_K$ such that
\be\label{measure}
      |\nu(\phi)| \; \le \; a^{}_K\|\phi\|
\ee
for all $\phi \in {\cal C}$ with support in $K$. Such measures are
in one-to-one correspondence with regular Borel measures through the
Riesz-Markov representation theorem, see \cite[Ch.\ XIII]{Dieu} and
\cite[Ch.\ 8, Sec.~69]{Berberian} for background material, 
and we thus identify these two pictures.
We deal only with regular Borel measures in this paper. 
If for each $\phi \in {\cal C}$ inequality (\ref{measure}) holds uniformly for 
all translates of $\phi$, we say that $\nu$ is {\em translation bounded}. 
This turns out to be a very useful concept because the {}Fourier transform of a 
tempered measure, though a tempered distribution, need not be a measure, but the 
Fourier transform of the autocorrelation of a translation bounded tempered 
distribution, if it exists, is always a measure. 
Let us explain how to control this subtle point in the context of diffraction theory.

Distributions and measures act on different spaces of functions and are equipped 
with different topologies. Nonetheless, there is an important connection between 
them which comes through the fact that the space ${\cal D}$ of $C^{\infty}$-functions 
of compact support is dense both in ${\cal C}$ and in ${\cal S}$. 

If a measure $\nu$ defines a tempered distribution $T_{\nu}$ by
\be
     (T_{\nu},\psi) \; = \; \nu(\psi) \; = \; \mbox{$\int \psi \, d\nu$}
\ee
for all $\psi\in\cal S$, the measure $\nu$ is called {\em tempered}.
A sufficient condition for a measure to be tempered is that it is
slowly increasing in the sense that 
$\int(1+|\bm x|)^{-k} |\nu|(d\bm x) < \infty$ for some $k\in\ZZ^+$,
see \cite[Thm.\ VII.VII]{Schwartz} or \cite[Ex.\ 7.12 b]{Rudin}.
Here, $|\nu|$ is the unique absolute value of $\nu$, i.e.\ the smallest positive
measure $\rho$ such that $|\nu(\phi)|\leq\rho(|\phi|)$ for all $\phi\in\cal C$.
It is also called the {\em variational} measure of $\nu$. 
As a partial converse, any positive tempered distribution
is a positive tempered measure. Thus, under the assumption of 
positiveness, tempered measures and tempered distributions can
be viewed as the same thing, and the Bochner-Schwartz theorem
can be restated as follows: 
$\nu$ is a positive tempered measure
if and only if it is the {}Fourier transform of a tempered distribution
of positive type.

Now, let $\nu$ be a translation bounded measure. Clearly $\nu$
is tempered. Our previous definition
of the autocorrelation of a Dirac comb has a natural extension to the
autocorrelation $\gamma^{}_{\nu}$ of $\nu$ by means of a volume-normalized
convolution of $\nu$ with $\nu^*$, where $\nu^*$ is defined by
$\nu^*(\phi) := \overline{\nu(\phi^*)}$, see \cite{Hof} for 
details\footnote{Note, however, that the definition of the analogous
operation to $*$ in \cite{Hof} is slightly incorrect in that it omits the
extra complex conjugation.}.
This construction of $\gamma^{}_{\nu}$ guarantees that, when it exists,
it is a positive definite tempered measure (and distribution), so combining the 
previous arguments we see that the corresponding diffraction pattern 
$\hat{\gamma}^{}_{\nu}$ is a tempered positive measure. 

This is important because it is this object that describes what one
actually sees on the screen in a diffraction experiment: each (measurable)
volume is assigned a non-negative number, namely the total intensity of
radiation scattered into this volume. In this article,
we will only meet the simple case that $\nu$ is a Dirac comb $\omega^{}_S$. 
However, already the convolution with a function of compact support (a ``profile''
of the scatterer) shows why the more general setting is useful. Let us
summarize this by the following result, which is a combination of
\cite[Prop.\ 2.2]{Hof} and \cite[Prop.\ 3.3]{Hof}.
\begin{prop} \label{tempmeasures} $\;$
  Let $\nu$ be a translation bounded measure. If its natural
  autocorrelation $\gamma^{}_{\nu}$ exists, it is a translation
  bounded (hence tempered) positive definite measure. {}Furthermore,
  $\hat{\gamma}^{}_{\nu}$ is then a positive measure and also
  translation bounded.
\end{prop}

\subsubsection*{Decomposition of measures}

A {\em pure point} of a measure $\nu$ is a point $\bm x\in\RR^n$ with
$\nu(\{\bm x\})\ne0$. Since $\nu$ is a regular Borel measure, it has
at most countably many pure points and the sum of $|\nu(\{\bm x\})|$ 
over the pure points in any compact set $K\/$ is convergent, and the pure points
alone serve to define a measure $\nu_{pp}$ called the {\em pure point part} of 
$\nu$. Thus $\nu$ has a unique decomposition as
\be \label{decomposition}
     \nu \; = \; \nu^{}_{pp}+\nu^{}_c\,,
\ee
where $\nu^{}_c := \nu - \nu^{}_{pp}$ is the so-called 
{\em continuous part}\footnote{The word ``continuous" here does not refer to being 
continuous as a function, since a line measure in the plane, for example, is 
continuous. It refers to the intermediate value property that if there are sets 
$A \subset C\/$ with $\nu_c(A)<b<\nu_c(C)$ then there is a set $B\/$ with 
$A \subset B \subset C\/$ and $\nu_c(B)=b$.
Some authors use the word diffuse instead of continuous and the words
atomic or purely discrete instead of pure point.}
of $\nu$, and is characterized by having no pure points.  
A {\em pure point measure} is a measure whose continuous part is 0. 

When a tempered measure $\nu$ is decomposed in this way 
the components $\nu_{pp}$ and $\nu_c$ are, of course, measures, but not necessarily
tempered.
{}For a translation bounded measure $\nu$, however, the components are both
translation bounded and hence tempered. 
If $\nu$ is a {\em positive} tempered measure, the decomposition is automatically
into tempered components.

Returning now to diffraction patterns, we meet the special situation
that $\hat{\gamma}^{}_{\omega}$ is a {\em positive} measure. Consequently,
it decomposes into a pure point part and a continuous part, both of which
are positive measures. The pure point part is called the {\em Bragg spectrum}
(of the point set $S$ that created it). We say that the diffraction pattern is
{\em pure point} or that $S$ has a {\em pure point diffraction spectrum} if its 
diffraction pattern is a pure point measure, i.e.\ if it is equal to
its Bragg spectrum.

Let us finally mention that, relative to Lebesgue measure, the continuous
part $\nu^{}_c$ can be further decomposed into an {\em absolutely continuous}
and a {\em singular continuous} part, $\nu^{}_c = \nu^{}_{ac} + \nu^{}_{sc}$,
see \cite{Berberian,RS} for details. In our examples, we show that $\nu_c$
vanishes, which means that there is neither an absolutely continuous nor
a singular continuous component present.

\subsection*{An intuitive derivation of the point spectra}

This section describes a short, intuitive way of calculating
the pure point part of the diffraction spectrum of the visible points $V\/$
of a lattice $\Gamma$ and also of the $k\/$th-power-free numbers $F=F_k$.
Its purpose is to give a taste of our later number theoretic methods in a
simpler setting and also to contrast this intuitively clear calculation
with the more circuitous route via autocorrelations we take later in rigorously 
establishing the complete diffraction spectra. It seems almost miraculous when 
the longer method eventually reduces to the same simple result.
The intuitive method depends on the fact that since the autocorrelation
$\gamma^{}_\omega$ is a volume-normalized convolution of
$\omega = \omega^{}_V$ with itself, its {}Fourier transform, the diffraction
pattern $\hat{\gamma}^{}_{\omega}$, should be a normalized square
of $\hat\omega$.  

\begin{figure}\centering
\setlength{\unitlength}{1mm}
\begin{picture}(40,28)
\put(0,0){\makebox(0,0){$\hat\omega$}}
\put(0,28){\makebox(0,0){$\omega$}}
\put(40,0){\makebox(0,0){$\hat{\gamma}^{}_\omega$}}
\put(40,28){\makebox(0,0){$\gamma^{}_\omega$}}
\put(20,4){\makebox(0,0){$|\cdot|^2$}}
\put(20,31){\makebox(0,0){$*$}}
\put(-5,14){\makebox(0,0){\sc ft}}
\put(45,14){\makebox(0,0){\sc ft}}
\put(4,0){\vector(1,0){31}}
\put(4,28){\vector(1,0){31}}
\put(0,25){\vector(0,-1){22}}
\put(40,25){\vector(0,-1){22}}
\end{picture}
\caption{Wiener diagram} \label{wiener}
\end{figure}
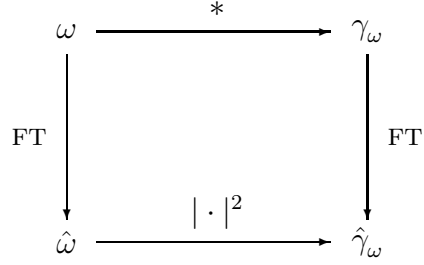

In this context, the coefficients of $\hat\omega$ are usually called
{\em amplitudes}, even if they only exist formally, while those
of $\hat{\gamma}^{}_{\omega}$ are called {\em intensities} which relates
to the fact that they are real and non-negative.
The appropriate operation now is to determine the
intensities as the absolute squares of the corresponding amplitudes
of peaks, as indicated in the commutative Wiener diagram in {}Figure~\ref{wiener}
where the vertical arrows represent the {}Fourier transform, the
upper horizontal arrow the volume-normalized convolution, and the
lower horizontal arrow taking the absolute squares of amplitudes.

This observation dates back at least 100 years in optics and is the
standard procedure in diffraction theory \cite{Cowley}. In the context
of diffraction from infinite arrangements, it
has been made rigorous (at least for the pure point part of the spectrum)
by Hof \cite{Hof,Hof2}. The intuitive approach to the visible
points runs through the Wiener
diagram along the ``low road", via the difficult to interpret
${\hat \omega}$.

\subsubsection*{Visible points}

We start from
\be\label{dist_sum}
  \omega^{}_V=\sum_{m=1}^\infty\mu(m)\omega^{}_{m\Gamma\setminus\{\bms 0\}} \, ,
\ee
which follows from (\ref{chiv}) since the sum on the right is supported
on $\Gamma\setminus\{\bm 0\}$ and the amplitude of the peak at
$\bm x\in\Gamma\setminus\{\bm 0\}$ is
\be\label{point_sum}
  \chi^{}_V(\bm x) \; = \; \sum_{m\mid\cont(\bms x)}\mu(m) \, .
\ee
Although the sum on the right of (\ref{point_sum}) is finite for each $\bm x$,
these pointwise sums are not uniformly convergent, so the sum on the right
of (\ref{dist_sum}) does not pass the $M\/$-test of Proposition~\ref{Mtest}
as a sum of pure point distributions.  Nevertheless, it is clear that the sum
does converge to $\omega^{}_V$ in the weak*-topology (though not in the local norm
topology).  Taking the {}Fourier transform of (\ref{dist_sum}) term-by-term we
obtain
\be\label{visibleFT}
   \hat\omega^{}_V \; = \;
   \sum_{m=1}^\infty\mu(m)\left({1\over m^n}\omega^{}_{\Gamma^*/m}-1\right),
\ee
where here and in the remainder of this subsection we assume $\Gamma$ to have
unit density, i.e.\ $\vol(\Gamma)=1$.
Since $\omega^{}_V$ is a tempered distribution it has a {}Fourier transform
$\hat\omega^{}_V$ which is also a tempered distribution, and since the
{}Fourier transform is a continuous operator on the space of tempered
distributions with the weak*-topology \cite{RS} the sum on the right of
(\ref{visibleFT}) converges to $\hat\omega^{}_V$ in the weak*-topology (though
again it does not pass the $M\/$-test).  This time the sums of the amplitudes
at individual points of $\QQ\Gamma^*$ are infinite sums but uniformly convergent,
when $n\ge2$.

It is difficult to interpret (\ref{visibleFT}) as a distribution and we suspect
that it is not a measure. Nevertheless, we can identify the pure points and
their amplitudes formally (even though they may not be real peaks).

Define the {\em denominator\/} of a non-zero point $\bm x\in\QQ\Gamma^*$ by
\be
      {\den}(\bm x) \; := \; \gcd\{m\in\ZZ\mid m\bm x\in \Gamma^*\}
\ee
(it is the smallest positive integer such that $m\bm x\in\Gamma^*$).
Then the sum of the amplitudes in (\ref{visibleFT}) at a point
$\bm x\in\QQ\Gamma^* \backslash \{\bm{0}\}$ with denominator $d\/$ is
\begin{eqnarray}
   \sum_{l=1}^\infty{\mu(ld)\over(ld)^n}
   & = & {\mu(d)\over d^n}\sum_\Cond{l=1\CR(l,d)=1}^\infty{\mu(l)\over l^n}
   \; = \; {\mu(d)\over d^n}\prod_{p\!\not\,\mid d}\left(1-{1\over p^n}\right)
   \nonumber\\
   &=&{\mu(d)\over\zeta(n)d^n}\prod_{p\mid d}\left(1-{1\over p^n}\right)^{-1}
   \; = \; {\mu(d)\over\zeta(n)}\prod_{p\mid d}{1\over p^n-1} \, .
\end{eqnarray}
When $n\ge2$, the pointwise sum of the pure point parts is thus nonzero at all
points of $\QQ\Gamma^*$ with squarefree denominator but is not absolutely
locally summable.  It is, however, locally square summable.  The squaring
operation gives, for the diffraction spectrum $\hat{\gamma}^{}_{\omega}$ of $V$,
the distribution with a pure point at each point of $\QQ\Gamma^*$ with squarefree
denominator, the peak at such a point with denominator $d\/$ having intensity
\be\label{intensity}
     {1\over\zeta^2(n)}\prod_{p\mid d}{1\over(p^n-1)^2} \, .
\ee
These coefficients are uniformly absolutely locally summable
and in fact are the correct intensities of the diffraction pattern. 
This will be proved below in Theorem \ref{diff_thm}.
They were first derived, based on similar arguments, in \cite{BGW}. 
However, one may well ask why this works.
How can we justify the squaring operation as the appropriate
mechanism for obtaining the intensities? The intuition is supported as follows.

\subsubsection*{Hof's results}

Two results of Hof \cite{Hof,Hof2} are
\begin{prop} $\;$
 Let $\nu$ be a translation bounded measure on $\RR^n$ such that $\hat{\nu}$ 
 is also a translation bounded measure. Then
\be
   \hat{\nu} (\{\bm x\}) \; = \; \lim_{R\to\infty}
        \frac{1}{v^{}_n R^n} \,
        \int_{B_R(\bms a)} e^{-2\pi i \bms x \cdot \bms y} \; \nu(d\bm y)
\ee
for every $\bm x\in\RR^n$ and the limit exists uniformly in $\bm a$.
\end{prop}

\begin{prop}\label{Hof2} $\;$
 Let $\nu$ be a translation bounded measure with natural autocorrelation $\gamma$
 and suppose that for all $\bm x \in \RR^n$,
\be\label{FTamp}
       m^{}_{\bms x} \; := \; \lim_{R\to\infty}
                              \frac{1}{v^{}_n R^n} \;
        \int_{B_R(\bms a)}  e^{-2\pi i \bms x \cdot \bms y}
                                    \; \nu(d\bm y) \, ,                     
\ee
  exists uniformly in $\bm a$. Then, for all $\bm x$, we have
\be
      \hat{\gamma}(\{\bm x\}) \; = \; |m^{}_{\bms x}|^2 \, .
\ee
\end{prop}
Taken together, these show that if the {}Fourier transform $\hat\nu$ of a
translation bounded tempered distribution $\nu$ is also translation bounded
(and hence can be decomposed as a pure point part and a continuous part) then
the pure points of the diffraction spectrum $\hat{\gamma}^{}_\nu$ are the same as
the pure points of the {}Fourier transform $\hat{\nu}$, but their intensities are 
the absolute squares of the amplitudes of $\hat{\nu}$.  Proposition~\ref{Hof2}
alone says that this continues to be true even when $\hat\nu$ is not translation
bounded, at least as long as the formal expressions for the pure point amplitudes,
also called {}Fourier-Bohr coefficients, of $\hat\nu$ 
(which may now not represent peaks of $\hat\nu$ in the accepted sense) 
converge uniformly with respect to translation of physical space.
We shall say more about the status of the method after giving 
the analogous intuitive derivation for the $k$th-power-free integers.

\subsubsection*{$k\/$th-power-free integers}

The parallel intuitive calculation for the pure point part of the 
diffraction spectrum of
the set $F=F_k$ of $k\/$th-power-free integers in $\RR$ goes like this.  We have
\be
 \omega^{}_F \; = \; 
   \sum_{m=1}^\infty\mu(m)\omega^{}_{m^k\ZZ\setminus\{0\}},
\ee
hence
\be\label{kfreeFT}
  \hat\omega^{}_F \; = \;
  \sum_{m=1}^\infty\mu(m)\left({1\over m^k}\omega^{}_{\ZZ/m^k}-1\right).
\ee
Let $d_k$ be the smallest positive integer such that $d_k^k$ is divisible by $d$.  
We note that $d_k$ is a divisor of $d$ but is divisible by every
prime factor of $d$ and that $d\mid m^k$ if and only if $d_k\mid m$.
Then the formal sum of the amplitudes of the peaks
of (\ref{kfreeFT}) at a point $x\in\QQ$ with denominator $d\/$ is
\be
  \sum_{l=1}^\infty{\mu(ld_k)\over(ld_k)^k}
  \; = \; {\mu(d_k)\over d_k^k}\sum_\Cond{l=1\CR(l,d)=1}^\infty{\mu(l)\over l^k}
  \; = \; {\mu(d_k)\over\zeta(k)d_k^k}\prod_{p\mid d}\left(1-{1\over p^k}\right)^{-1}
  \; = \; {\mu(d^*)\over\zeta(k)}\prod_{p\mid d}{1\over p^k-1}\, ,
\ee
since, when $d_k$ is squarefree, we have $d_k=d^*=\prod_{p\mid d}p\/$ (the
squarefree kernel of $d\/$) and $d_k^k=\prod_{p\mid d}p^k$.  Consequently,
the pure point part of the diffraction spectrum of $F\/$ has peaks at the
points of $\QQ$ with $(k+1)$th-power-free denominator, and its intensity
at a such a point, with denominator $d\/$, is
\be
{1\over\zeta^2(k)}\prod_{p\mid d}{1\over(p^k-1)^2}.
\ee
This is a pure point distribution and also agrees with the result of
Theorem~\ref{diff_thm_F} below.

\subsubsection*{Status of the method}

In the remaining sections of this paper, which contain the heart
of the proofs, we completely avoid Hof's formula (\ref{FTamp}) 
(which in any case could only determine the pure point part of 
the spectrum without any assurance that there is no continuous part).  
What prevents us from making our derivation of the pure point
part of the spectrum rigorous by citing Proposition~\ref{Hof2} is that for
our examples the limit in (\ref{FTamp}) is not uniform in \bm a.  This is
most easily seen by noting that our point sets have arbitrarily large holes,
so that no matter how large $R\/$ is there will be an \bm a such that
$B_R(\bm 0)+\bm a$ lies entirely within a hole.  However, large holes are
very sparse, so the limit is, in a sense, ``nearly uniform".  {}For 
$\bm x \not\in \QQ\Gamma^*$
the limit (\ref{FTamp}) can be shown to be 0 (as it should be) but
again is almost certainly not uniform in \bm a.

The fact that our examples are outside the domain where Hof's result applies
underscores the fact that we are in new territory not only for calculating
the continuous part of the diffraction spectrum but even for calculating the
diffraction peaks. In other words: the answer to the question which distributions
of matter diffract is still largely unknown.

\subsection*{Autocorrelation of the visible points}

The remainder of this paper consists of a rigorous derivation of the results
of the previous section via the ``high road" of the autocorrelation.  So we
begin by using an elaboration of the kind of M\"obius inversion argument
used in proving Proposition~\ref{freq1} to calculate the autocorrelation
of the visible points $V\/$ of a lattice $\Gamma$.

\begin{theorem}\label{ACthm} $\;$
{}For $n\ge3$, the natural autocorrelation $\gamma$ of the set $V\/$ 
of visible points
of a lattice $\Gamma$ exists and is supported on $\Gamma$, the weight of
a point $\bm a\in\Gamma$ in the autocorrelation of $V\/$ being given by
\be  \label{thm1-eq}
   w(\bm a) \; = \; {\rm dens}(\Gamma)\,\xi(n)\!\!\!
        \prod_{p\mid\cont(\bms a)}\left(1+{1\over p^n-2}\right),
\ee
with error term (as defined just after (\ref{density})) equal to $O(1/R)$,
where the implied constant depends on \bm a as well as on $\Gamma$.
\end{theorem}
{\sc Proof}: By Prop.~\ref{diffprop}, $V-V=\Gamma$, so the
autocorrelation of $V\/$ (if it exists) is supported on $\Gamma$.
The weight of a point $\bm a\in\Gamma$ in the autocorrelation
of $V\/$ is the limit as $R\to\infty$ of
\be \label{sum}
{1\over v_nR^n}\!\!
\sum_\cond{|\bms x|,|\bms x-\bms a|< R\CR\bms x,\bms x-\bms a\in V}
\!\!\!\!\!\!1
\ee
and by Lemma~\ref{ACexist} the existence of this limit for each
$\bm a\in\Gamma$ is sufficient to ensure the existence of the autocorrelation.

It is convenient to drop the condition $|\bm x-\bm a|<R\/$ in (\ref{sum}),
which then becomes
\be\label{ACvis}
  {1\over v_nR^n}\;\sum_\Cond{|\bms x|< R\CR\bms x,\bms x-\bms a\in V}\;1 \, .
\ee
The difference between these sums is $O(1/R)$, due to the extra lattice points
\bm x within a constant distance $|\bm a|$ of the boundary of $B_R(\bm0)$ that
are included in the latter.  We note that the latter sum has a natural geometric
interpretation as the proportion of integer points in a large ball that are
visible both from the origin and from the viewpoint \bm a.  By (\ref{chiv}) and
(\ref{contbound}), (\ref{ACvis}) can be expressed as
\be\label{sumACvis}
  {1\over v_nR^n}\!\!
  \sum_\cond{\bms x\in\Gamma\setminus\{\bms0,\bms a\}\CR|\bms x|<R}
  \sum_\cond{l\mid\cont{(\bms x)}}\mu(l)\sum_\Cond{m\mid\cont{(\bms x-\bms a)}}
  \mu(m)
  \; = \; {1\over v_nR^n}\!\!\sum_\cond{1\le l<S}\sum_\cond{1\le m<S}\mu(l)\mu(m)
  \sum_\Cond{\bms x\in\Gamma\setminus\{\bms0,\bms a\}\CR
  \bms x\equiv\bms0\mod{l\Gamma}\CR\bms x\equiv\bms a\mod{m\Gamma}\CR
  |\bms x|<R}\;1 \, ,
\ee
where $S=(R+|\bm a|)/L(\Gamma)$.  Collecting together terms with the same
value of $d=(l,m)$, noting that all solutions \bm x of the congruences
belong to $d\Gamma$ and that the congruences have no solution at all unless
$\bm a\in d\Gamma$, and putting $l'=l/d$, $m'=m/d$, $\bm x'=\bm x/d$,
$\bm a'=\bm a/d$, we obtain
\be\label{triplesum}
  {1\over v_nR^n}\!\!\sum_\cond{d\mid\cont(\bms a)}\sum_\cond{1\le l'<S/d}
  \sum_\cond{1\le m'<S/d\CR(l',m')=1}\mu(l'd)\mu(m'd)
  \sum_\Cond{\bms x'\in\Gamma\setminus\{\bms0,\bms a'\}\CR
  \bms x'\equiv\bms0\mod{l'\Gamma}\CR\bms x'\equiv\bms a'\mod{m'\Gamma}\CR
  |\bms x'|<R/d}\;1 \, .
\ee
Since $l'$ and $m'$ are bound variables of summation and $\bm x'$ and $\bm a'$
will not be referred to again, we can drop the dashes: from now on $l\/$ and
$m\/$ are the new $l'$ and $m'$ but \bm a is the original \bm a.

By Propositions~\ref{CRT} and~\ref{sphere} (with $lm\Gamma$ in place of $\Gamma$)
the inmost sum is
\be
   {\rm dens}(\Gamma)\,v_n\left({R\over dlm}\right)^n
   +O\left({R\over dlm}\right)^{n-1}+O(1)
\ee
These three terms give a main term and two error terms in (\ref{triplesum}).

The first error term is majorized by
\be
  O\left({1\over R}
  \sum_{1\le l<S}{1\over l^{(n-1)}}\sum_{1\le m<S}{1\over m^{(n-1)}}\right)
  \; = \; \cases{O(1/R),&if $n\ge3$,\cr O((\log S)^2/R),&if $n=2$,\cr}
\ee
since the sums are convergent when $n-1\ge2$ and increase logarithmically
when $n-1=1$.

The second error term is majorized by $O(S^2/R^n)=O(1/R^{n-2})$.
So when $n\ge3$ both error terms are $O(1/R)$ (since $S=O(R)$) and
tend to 0 as $R\to\infty$.

The main term is (with ${\cal D} = {\rm dens}(\Gamma)$)
\be
  {\cal D}\!\!\!\!\!\sum_\cond{d\mid\cont(\bms a)}\!\sum_\cond{1\le l<S/d}
  \!\sum_\cond{1\le m<S/d\CR(l,m)=1}\!\!\!\!{\mu(ld)\mu(md)\over(dlm)^n}
  \; = \; {\cal D}\!\!\!\!\!
  \sum_\cond{d\mid\cont(\bms a)}\!\sum_\cond{1\le l<S/d\CR(l,d)=1}
  \!\sum_\cond{1\le m<S/d\CR(m,d)=1\CR(l,m)=1}\!\!\!\!
  {\mu^2(d)\mu(l)\mu(m)\over(dlm)^n}
\ee
since $\mu(ld)$ is $\mu(l)\mu(d)$ when $(l,d)=1$ and 0 otherwise,
and similarly for $\mu(md)$,
\be\label{partial}
  = \; {\cal D}\!\!\!\!\sum_\cond{d\mid\cont(\bms a)}{\mu^2(d)\over d^n}
  \sum_\cond{1\le l<S/d\CR(l,d)=1}
  \sum_\cond{1\le m<S/d\CR(m,d)=1}{\mu(lm)\over(lm)^n}
\ee
since $\mu(l)\mu(m)$ is $\mu(lm)$ when $(l,m)=1$ and 0 otherwise
\be\label{limit}
  \longrightarrow \;\; 
  {\cal D}\!\!\!\!\sum_\cond{d\mid\cont(\bms a)}{\mu^2(d)\over d^n}
  \sum_\Cond{r=1\CR(r,d)=1}^\infty{\mu(r)\sigma(r)\over r^n}
  \quad\quad\mbox{ as }R\to\infty
\ee
since the double sum is absolutely convergent.  The difference between
this limit and the partial sum (\ref{partial}) is $O(1/R^{n-1})$, so
falls within the error term estimate $O(1/R)$.

Using (\ref{xi}), the expression for the limit (\ref{limit}) can be rearranged as
\begin{eqnarray}
  \lefteqn{
  {\cal D}\;\sum_\Cond{d\mid\cont(\bms a)\CR d{\rm\ squarefree}}\;\;\;
  {1\over d^n}
  \prod_{p\!\not\;\mid d}\left(1-{2\over p^n}\right)
  \;=\;{\cal D}\,\xi(n)
  \sum_\Cond{d\mid\cont(\bms a)\CR d{\rm\ squarefree}}\;\;\;{1\over d^n}
  \prod_{p\mid d}\left(1-{2\over p^n}\right)^{-1}}\label{AC}\\
  &=&{\cal D}\,\xi(n)\!\!\!\!\prod_{p\mid\cont(\bms a)}
  \left(1+{1\over p^n}\left(1-{2\over p^n}\right)^{-1}\right)
  \;=\;{\cal D}\,\xi(n)\!\!\!\!
  \prod_{p\mid\cont(\bms a)}\left(1+{1\over p^n-2}\right). \nonumber
\end{eqnarray}
This completes the proof. \hfill $\square$

To establish the existence of the autocorrelation function for the
visible points when $n=2$, we need to reduce the second error term in
the above proof. This we do by modifying the argument slightly, first
using the characteristic function (\ref{chiv}) to replace the constraint
on $\bm x-\bm a$ only, then discarding large values of $m\/$ from the sum
before using (\ref{chiv}) again to replace the constraint on~$\bm x$.

\begin{theorem}\label{ACthm2} $\;$
   Theorem~\ref{ACthm} holds for $n=2$ with the error term increased to
   $O(1/\sqrt R)$.
\end{theorem}
{\sc Proof}: Using (\ref{chiv}) to replace the constraint on $\bm x-\bm a$,
(\ref{ACvis}) becomes
\be\label{vis2}
  {1\over\pi R^2}\!\!\sum_\cond{\bms x\in V\CR\bms x\ne\bms a\CR|\bms x|<R}
  \sum_\cond{m\mid\cont(\bms x-\bms a)}\mu(m)
  \; = \; {1\over\pi R^2}\!\!\sum_\cond{1\le m<S}\mu(m)
  \sum_\Cond{\bms x\in V\CR\bms x\ne\bms a\CR
  \bms x\equiv\bms a\mod{m\Gamma}\CR|\bms x|<R}\;\;1 \, ,
\ee
where $S\/$ is as before.  The inner sum is trivially $O(R^2/m^2)$
(since $m<S=O(R)$) so the contribution to (\ref{ACvis}) from the terms
with $m\ge\sqrt R$ is
\be
  O\biggl(\sum_\cond{m\ge \sqrt R}{1\over m^2}\biggr)\;=\;O(1/\sqrt R) \, .
\ee
{}For the terms in (\ref{vis2}) with $m<\sqrt R\/$ we use the characteristic
function (\ref{chiv}) to replace the constraint on $\bm x$, obtaining
\be
  {1\over\pi R^2}\!\!\sum_\cond{1\le m<\sqrt R}\mu(m)
  \sum_\Cond{\bms x\in\Gamma\setminus\{\bms 0,\bms a\}\CR
  \bms x\equiv\bms a\mod{m\Gamma}\CR|\bms x|<R}\qquad
  \sum_\Cond{l\mid\cont(\bms x)}\mu(l) \; = \;
  {1\over\pi R^2}\!\!\sum_\cond{1\le l<S}\sum_\cond{1\le m<\sqrt R}\mu(l)\mu(m)
  \sum_\Cond{\bms x\in\Gamma\setminus\{\bms 0,\bms a\}\CR
  \bms x\equiv\bms 0\mod{l\Gamma}\CR\bms x\equiv\bms a\mod{m\Gamma}\CR|\bms x|<R}1\,.
\ee
This is now identical to (\ref{sumACvis}), except that the second $S\/$
is replaced by $\sqrt R\/$, and as before contributes to (\ref{ACvis})
a main term that tends to
\be
   {\cal D}\,\xi(2)\!\!\prod_{p\mid\cont(\bms a)}\left(1+{1\over p^2-2}\right)
\ee
as $R\to\infty$ and two error terms $O((\log R)^2/R)$ and $O(1/\sqrt R)$.
The difference between the main term and its limit is also $O(1/\sqrt R)$.
So the total error term is majorized by $O(1/\sqrt R)$. \hfill $\square$

Note that this modified argument does not reduce the error term
when $n\ge3$, which contains a term $O(1/R)$ arising from
the boundary of the ball when $l=m=1$.
{}For $n=2$, a similar but more complicated argument gives an improved
error term $O(R^{-3/4}(\log R)^c)$ for some constant $c$, but we do not need
this here.

{\bf Remark}\quad Theorems~\ref{ACthm} and \ref{ACthm2} show that the weight
$w(\bm a)$ in the autocorrelation of $V\/$ (also called the autocorrelation
coefficient\footnote{Note that $w(\bm a)$ can also be interpreted as the density
of lattice points that are simultaneously visible both from the origin and from
the lattice point $\bm a$. 
Extending this to the condition of simultaneous visibility from an arbitrary
(but finite) set of points of $\Gamma$ results in higher order correlation
coefficients, see \cite{Berend} for some recent results.}) depends only on the 
content of $\bm a$ and the density of $\Gamma$. 
This enables us to calculate how the application of a non-singular linear
transformation  $T\/$ to $V\/$ affects the autocorrelation of $V$.
Clearly $T\/$ preserves content, in the sense that, for any lattice
$\Gamma$ and $\bm x\in\Gamma$, the content of $T\bm x$ as a vector of
$T\Gamma$ is equal to the content of $\bm x$ as a vector of $\Gamma$.
In particular, $TV\/$ is the set of visible points of $T\Gamma$. The 
autocorrelation of $TV\/$ is supported on $T\Gamma$ and (\ref{thm1-eq}) shows that
\be\label{ACtransform}
  w^{}_{TV}(T\bm a) \; = \; {\dens(T\Gamma)\over\dens(\Gamma)}\;w^{}_V(\bm a)
   \; = \; {1\over|\det T|}\;w^{}_V(\bm a) \, ,
\ee
where the suffix on $w\/$ indicates which autocorrelation it is associated
with.  Thus the autocorrelation of $TV\/$ is $|\det T|^{-1}$ times the
$T\/$-image of the autocorrelation of $V$.

The direct way of calculating $w_{TV}(T\bm a)$ from (\ref{sum}) is
\be
     w^{}_{TV}(T\bm a) \; = \; \lim_{R\to\infty}{1\over v_nR^n}
       \sum_{\bms x,\bms x-T\bms a\in R B\cap TV}\!\!\!\!\!1 \, ,
\ee
where $B=B_1(\bm 0)$ is the unit ball in $\RR^n$, and substituting this in
(\ref{ACtransform}) gives
\be
  w^{}_V(\bm a) \; = \; \lim_{R\to\infty}{|\det T|\over v_nR^n}
  \sum_{\bms y,\bms y-\bms a\in R T^{-1}B\cap V}\!\!\!\!\!1 \; = \;
  \lim_{R\to\infty}{1\over\vol(E)R^n}
  \sum_{\bms y,\bms y-\bms a\in R E\cap V}\!\!\!\!\!1  \, ,
\ee
where $E=T^{-1}B\/$ is an ellipsoid.  Hence averaging over the expanding
ellipsoid $RE\/$ gives the same value for the autocorrelation of $V\/$ as
using the natural density and averaging over an expanding ball.  The
choice of ellipsoid here is completely arbitrary, since for any $E\/$
there is a linear transformation $T\/$ with $TE=B\/$.  In the Appendix
this is generalized further and we show how to replace $E\/$ by any
bounded measurable region (not necessarily centred at $\bm 0$) with
finite $(n-1)$-dimensional surface area.

\subsection*{Diffraction spectrum of the visible points}

The final step in obtaining the diffraction spectrum of the visible points
$V\/$ is to take the {}Fourier transform of the autocorrelation of $V$. 

\begin{theorem}\label{diff_thm} $\;$
   The diffraction spectrum of the set of visible points of an
   $n$-dimensional lattice $\Gamma$ (with $n\ge 2$)
   exists and is a pure point measure which is concentrated on the set
   of points in $\QQ\Gamma^*$ with squarefree denominator and whose intensity at
   a point with such a denominator $q\/$ is given by
\be\label{diff_amp}
   {\dens(\Gamma)^2\over\zeta^2(n)}\;\prod_{p\mid q}{1\over(p^n-1)^2} \, .
\ee
This measure can also be represented as
\be\label{diff_sum}
  {\rm dens}(\Gamma)^2\,\xi(n)\sum_\Cond{d=1\CR d{\rm\ squarefree}}^\infty\;
  \biggl(\prod_{p\mid d}{1\over p^{2n}-2p^n}\biggr)\omega^{}_{\Gamma^*/d} \, ,
\ee
a weak*-convergent sum of Dirac combs.
\end{theorem}
{\sc Proof\/}: Let $\gamma$ be the autocorrelation of $V$.  The right hand
side of (\ref{AC}) can be expressed in the form 
(again, with ${\cal D} = {\rm dens}(\Gamma)$)
\be
  w(\bm a) \;= \;{\cal D}\,\xi(n)
  \sum_\Cond{d=1\CR d{\rm\ squarefree}\CR\bms a\in d\Gamma}^\infty\;{1\over d^n}
  \prod_{p\mid d}\left(1-{2\over p^n}\right)^{-1}.
\ee
So by Theorems~\ref{ACthm} and~\ref{ACthm2} and Lemma~\ref{ACexist}
\be
  \gamma \; = \; {\cal D}\,\xi(n)
  \sum_\Cond{d=1\CR d{\rm\ squarefree}}^\infty\;{1\over d^n}
  \prod_{p\mid d}\left(1-{2\over p^n}\right)^{-1}\omega^{}_{d\Gamma} \, .
\ee

Since $\|\omega^{}_{d\Gamma}\|^{}_\loc=O(1)$ and the coefficient of
$\omega^{}_{d\Gamma}$ is $O(1/d^n)$ when $n\ge2$, this sum of tempered
distributions is convergent in the weak*-topology by Lemma~\ref{Mtest}
(and this is easy to see by other means in this case, where the resulting
sum is uniformly discrete).  Its term-by-term {}Fourier transform is
\be
  {\cal D}^2\,\xi(n)
  \sum_\Cond{d=1\CR d{\rm\ squarefree}}^\infty\;{1\over d^{2n}}
  \prod_{p\mid d}\left(1-{2\over p^n}\right)^{-1}\omega^{}_{\Gamma^*/d},
\ee
which weak*-converges to the diffraction spectrum of $V$, since
the {}Fourier transform operator is weak*-continuous.  Since
$\|\omega^{}_{\Gamma^*/d}\|^{}_\loc=O(d^n)$ and the coefficient of
$\omega^{}_{\Gamma^*/d}$ is $O(1/d^{2n})$, Lemma~\ref{Mtest} tells
us that the weak*-sum is a translation bounded pure point measure
equal to the pointwise sum of its terms.

This establishes the series form (\ref{diff_sum}) for the diffraction spectrum.
The explicit values of the intensities can now be evaluated
quite simply. Let \bm p be a point in $\QQ\Gamma^*$ with
denominator $q$.  
We can suppose that $q$ is square-free, since otherwise there
is no contribution to (\ref{diff_sum}) at all.  
The terms in (\ref{diff_sum}) that contribute to the
intensity at \bm p are those with $d=mq\/$ ($m\/$ squarefree in $\ZZ^+$ and
prime to $q\/$), so the intensity at \bm p is
\be\label{point_amp}
  {\cal D}^2\,\xi(n)\prod_{p\mid q}{1\over p^{2n}-2p^n}
  \sum_\Cond{m=1\CR m{\rm\ squarefree}\CR(m,q)=1}^\infty\;\;
  \prod_{p\mid m}{1\over p^{2n}-2p^n}.
\ee
This simplifies to
\begin{eqnarray}
  \lefteqn{{\cal D}^2\,\xi(n)\prod_{p\mid q}{1\over p^{2n}-2p^n}
  \prod_{p\!\not\,\mid q}\biggl(1+{1\over p^{2n}-2p^n}\biggr)}\nonumber\\
  &=&{\cal D}^2\,\xi(n)
  \prod_{p\mid q}{1\over p^{2n}}\biggl(1-{2\over p^n}\biggr)^{-1}
  \prod_{p\!\not\,\mid q}
  \biggl(1-{1\over p^n}\biggr)^2\biggl(1-{2\over p^n}\biggr)^{-1}\nonumber\\
  &=&{{\cal D}^2\over\zeta^2(n)}
  \prod_{p\mid q}{1\over p^{2n}}\biggl(1-{1\over p^n}\biggr)^{-2}
\end{eqnarray}
(using the Euler products in (\ref{invzeta}) and (\ref{xi})) which
agrees with (\ref{diff_amp}). \hfill $\square$

An explicit example of the diffraction (for $\Gamma=\ZZ^2$) is shown in 
\cite{BGW}, and compared with an optical experiment. The diffraction image
in this case is both $D_4$-symmetric and $GL(2,\ZZ)$ invariant which results
in a beautiful image with a rather unusual symmetry structure, reminiscent
of self-similar patterns common in fractals.

\subsection*{$k\/$th-power-free numbers}

In this section we derive the diffraction spectrum of the 1-dimensional
set $F\/$ consisting of the $k\/$th-power-free numbers in $\ZZ$.  Again,
this has arbitrarily long gaps but it nevertheless has a pure point
diffraction spectrum.  The proof of the second assertion closely parallels
the corresponding proof for the visible points $V\/$ in the previous two
sections, with the parameter $k\/$ for $F\/$ playing the r\^ole of the
dimension $n\/$ of $V\/$ in the formalism.  There are some differences
of detail, however, particularly with the error terms.

Let us note here that the results on the autocorrelation derived below also
follow from \cite{Pillai} ($k=2$) and from \cite[Thm.\ 1]{Mirsky}, where
also more general correlation functions have been derived. 
The error term given in \cite{Mirsky} is $O(R^{-1+2/(k+1)})$, which is
slightly better than the error term we derive in Theorem~\ref{ACthmkf}.
We include our proof to make this paper self-contained and to show up
the close parallel between visible lattice points and $k$-free numbers.

\begin{prop} $\;$
   $F\/$ is uniformly discrete, but has gaps of arbitrary length.  Moreover,
for any $L>0$ the set of gaps of length at least $L\/$ has positive density.
\end{prop}
{\sc Proof}: The uniform discreteness is trivial since $F\subset\ZZ$.  
Choose $L\/$ integer moduli $m^{}_1,\ldots,m^{}_L$ that are $>1$
and coprime in pairs (for example, the first $L\/$ primes).  By the Chinese
Remainder Theorem there is an integer $N\/$ with
\be
   N\equiv-j+1\mod{m_j^k}\quad\mbox{for }j=1,\ldots,L.
\ee
Now for $x\equiv N\mod{m_1^km_2^k\ldots m_L^k}$ we have $m_1^k\mid x$,
$m_2^k\mid(x+1)$,\dots, $m_L^k\mid(x+L-1)$, so none of the numbers
$x,x+1,\ldots,x+L-1$ is $k\/$th-power-free.  The numbers $x\/$ have density
$(m^{}_1 m^{}_2\ldots m^{}_L)^{-k}$. \hfill $\square$

This argument gives a distance of the order $L^{kL}$ between gaps of
length $L$, so again long gaps can be expected to be extremely sparse.
Nevertheless, gaps of length $L$ have a definite frequency
(its expression in terms of Dirichlet series can be extracted from \cite{Mirsky}).
It is interesting to note that the corresponding distribution is {\em not}
Poissonian, but that these frequencies decline faster than exponentially in $L$
\cite{Huxley,Grimmett}.

\begin{prop} \label{prop11} $\;$
  {}For $k\ge2$, the $k\/$th-power-free integers $F\/$ have a natural density given by
\be
        {\rm dens}(F) \; = \; \frac{1}{\zeta(k)} \, ,
\ee
  with error term $O(R^{-1+(1/k)})$.
\end{prop}
At least for the squarefree numbers, this is again a
standard example of M\"obius inversion, see \cite[Thm.~6.6.1]{Hua}.

{\sc Proof}:
The natural density of $F\/$ is the limit as $R\to\infty$ of
\be\label{Dkf}
     {1\over 2R}\sum_\Cond{|x|< R\CR x\in {}F} 1 \, ,
\ee
which by (\ref{chik}) is
\begin{eqnarray}
     {1\over 2R}\sum_\cond{|x|< R\CR x\ne0}
                \sum_\Cond{m^k\mid x}\mu(m)
  &=&{1\over 2R}\sum_\cond{1\le m<R^{1/k}}\mu(m)
          \sum_\Cond{|x|< R\CR x\equiv 0\mod {m^k}\CR x\ne0} \;1 \label{sumDkf}\\
  &=&\sum_\cond{1\le m<R^{1/k}}\!{\mu(m)\over m^k}
          +O\left({R^{1/k}\over R}\right),\nonumber
\end{eqnarray}
since the inner sum on the right of (\ref{sumDkf}) is $2R/m^k+O(1)$.
This last expression tends to $1/\zeta(k)$ as $R\to\infty$ when $n\ge2$,
with the errors due to the tail of the sum and to the explicit error
term both $O(R^{-1+(1/k)})$. \hfill $\square$

\begin{theorem}\label{ACthmkf} $\;$
   {}For $k\ge2$, the natural autocorrelation of the set $F\/$ of $k\/$th-power-free
   integers exists and is supported on $\ZZ$, the weight of an integer $a\/$
   in the autocorrelation of $V\/$ being given by
\be
   w(a)\;=\;\xi(k)\prod_{p^k\mid a}\left(1+{1\over p^k-2}\right),
\ee
with error term $O(R^{-(1-(1/k))^2})$, where the implied constant
depends on $a$.
\end{theorem}
{\sc Proof}: Clearly the autocorrelation is supported on $\ZZ$.  The weight of
an integer $a\/$ in the autocorrelation of $F\/$ is the limit as $R\to\infty$ of
\be
{1\over2R}\!\!
\sum_\cond{|x|,|x-a|< R\CR x,x-a\ k\mbox{\footnotesize -}{\rm free}}\!\!\!\!\!\!1
\ee
and, as in the proof of Theorem~\ref{ACthm}, to show the existence of the
autocorrelation it is enough to show that this limit exists for every $a$.
Again as in the proof of Theorem~\ref{ACthm}, we can drop the condition
$|x-a|<R\/$ from the sum to obtain
\be\label{ACkfn}
{1\over2R}\;\sum_\Cond{|x|< R\CR x,x-a\ k\mbox{\footnotesize -}{\rm free}}\;1
\ee
with error $O(1/R)$.

This can be evaluated by an argument exactly parallel to that for the visible
points, using the characteristic function (\ref{chik}) in place of
(\ref{chiv}).  We use the argument in the modified form we used in the proof
of Theorem~\ref{ACthm2} since, in the context of $k\/$th-power-free numbers,
it gives an improved error term in all cases owing to the fact that, for a
1-dimensional set, errors due to the boundary of a large region are trivial. 
(Using the unmodified form would give error term $O(R^{-1+(2/k)})$ instead of
the $O(R^{-1+(2/k)-(1/k)^2})$ we obtain here.  As with the visible points,
this would not be small enough to establish the existence of the limit for
the squarefree numbers $k=2$ --- the paradigm case.)

Using the characteristic function (\ref{chik}) to replace the constraint
on $x-a$, (\ref{ACkfn}) becomes
\be\label{sumkfn}
  {1\over2R}\sum_\cond{|x|<R\CR x\ k\mbox{\footnotesize -}{\rm free}\CR x\ne a}
  \sum_\cond{m^k\mid(x-a)}\mu(m)
  \;=\;{1\over2R}\sum_\cond{1\le m<S^{1/k}}\mu(m)
  \sum_\Cond{|x|<R\CR x\ k\mbox{\footnotesize -}{\rm free}
  \CR x\ne a\CR x\equiv a\mod{m^k}}\;1\, ,
\ee
where $S=R+|a|$.  The inner sum is trivially $O(R/m^k)$ (since $m^k<S\/$) so
the contribution to (\ref{ACkfn}) from the terms with $m\ge T$, where $T\/$
is a parameter tending to infinity with $R\/$ to be chosen later, is
\be
  O\biggl(\sum_\cond{m\ge T}{1\over m^k}\biggr)\;=\;O(T^{1-k})\, .
\ee
{}For the terms in (\ref{sumkfn}) with $m<T\/$ we use the characteristic
function (\ref{chik}) to replace the constraint on $x$, obtaining
\be
  {1\over2R}\sum_\cond{1\le m<T}\mu(m)
  \sum_\Cond{|x|<R\CR x\ne0,a\CR x\equiv a\mod{m^k}}\;\;\sum_\cond{l^k\mid x}\mu(l)
  \;=\;{1\over2R}\sum_\cond{1\le l<R^{1/k}}\sum_\cond{1\le m<T}\mu(l)\mu(m)
\sum_\Cond{|x|<R\CR x\ne0,a\CR x\equiv 0\mod{l^k}\CR x\equiv a\mod{m^k}}\;1\, .
\ee
Collecting together terms with the same value of $d=(l,m)$, noting
that all solutions $x\/$ of the congruences are divisible by $d^k$
and that the congruences have no solution unless $d^k\mid a$,
and putting $l'=l/d$, $m'=m/d$, $x'=x/d^k$, $a'=a/d^k$, we obtain
\be
  {1\over2R}\sum_\cond{d^k\mid a}\sum_\cond{1\le l'<R^{1/k}/d}
  \sum_\cond{1\le m'<T/d\CR(l',m')=1}\mu(l'd)\mu(m'd)
  \sum_\Cond{|x'|< R/d^k\CR x'\ne0,a'\CR
  x'\equiv0\mod{l'^k}\CR x'\equiv a'\mod{m'^k}}\;1\, .
\ee
As in the proof of Theorem~\ref{ACthm}, we can drop the dashes from now on.

By the Chinese Remainder Theorem the inner sum is
\be
    {2R\over(dlm)^k}+O(1)\, ,
\ee
giving a main term and another error term in (\ref{ACkfn}).

The error term is majorized by $O(R^{-1+(1/k)}T)$.

The main term is
\begin{eqnarray}
  \lefteqn{\sum_\cond{d^k\mid a}\sum_\cond{1\le l<R^{1/k}/d}
  \sum_\cond{1\le m<T/d\CR(l,m)=1}{\mu(ld)\mu(md)\over(dlm)^k}}\nonumber\\
  &=&{\sum_\cond{d^k\mid a}\sum_\cond{1\le l<R^{1/k}/d\CR(l,d)=1}
  \sum_\cond{1\le m<T/d\CR(m,d)=1\CR(l,m)=1}{\mu^2(d)\mu(l)\mu(m)\over(dlm)^k}}
\end{eqnarray}
since $\mu(ld)$ is $\mu(l)\mu(d)$ when $(l,d)=1$ and 0 otherwise,
and similarly for $\mu(md)$,
\be\label{partial_kf}
  =\;\sum_\cond{d^k\mid a}{\mu^2(d)\over d^k}
  \sum_\cond{1\le l<R^{1/k}/d\CR(l,d)=1}
  \sum_\cond{1\le m<T/d\CR(m,d)=1}{\mu(lm)\over(lm)^k}
\ee
since $\mu(l)\mu(m)$ is $\mu(lm)$ when $(l,m)=1$ and 0 otherwise
\be
  \longrightarrow\;\;\sum_\cond{d^k\mid a}{\mu^2(d)\over d^k}
  \sum_\Cond{r=1\CR(r,d)=1}^\infty{\mu(r)\sigma(r)\over r^k}
  \quad\quad\mbox{ as }R\to\infty
\ee
since the double sum is absolutely convergent.  The difference
between this limit and the partial sum (\ref{partial_kf}) is
$O(1/R^{-(k-1)/k}) +O(T^{-(k-1)})$ and each of these terms is
subsumed by one of the error estimates we already have.

This last expression is
\begin{eqnarray}
  \lefteqn{\sum_\Cond{d^k\mid a\CR d{\rm\ squarefree}}{1\over d^k}
  \prod_{p\!\not\,\mid d}\left(1-{2\over p^k}\right)
  \;=\;\xi(k)\sum_\Cond{d^k\mid a\CR d{\rm\ squarefree}}{1\over d^k}
  \prod_{p\mid d}\left(1-{2\over p^k}\right)^{-1}}\label{ACkf}\\
  &=&\xi(k)\prod_{p^k\mid a}
  \left(1+{1\over p^k}\left(1-{2\over p^k}\right)^{-1}\right)
  \;=\;\,\xi(k)\prod_{p^k\mid a}\left(1+{1\over p^k-2}\right).\nonumber
\end{eqnarray}
{}Finally, we choose $T=R^{(1/k)-(1/k)^2}$, making both error terms
$O(R^{-1+(2/k)-(1/k)^2})$. \hfill $\square$

We can now fill a small gap in our analogy by the following result.
\begin{coro}
  $\quad$ {}For $k\ge2$, we have $F-F=\ZZ$.
\end{coro}
{\sc Proof\/}: Since $F-F\subset\ZZ$, $\gamma^{}_{\omega}$ is clearly
supported on $\ZZ$. {}From Theorem~\ref{ACthmkf}, we get $w(a)>0$ for
all $a\in\ZZ$, hence also $\ZZ\subset {}F-F$. \hfill $\square$

Let us finally describe the diffraction of $k\/$th-power-free integers.
\begin{theorem}\label{diff_thm_F} $\;$
 The diffraction spectrum of the set $F\/$ of $k\/$th-power-free integers
 exists and is a pure point measure. It is supported on the
 set of numbers $a/q\in\QQ$ with $q\/$ being $(k+1)$-power-free.
 The diffraction intensity at a point with such a denominator $q\/$ is
\be\label{diff_amp_F}
   {1\over\zeta^2(k)}\prod_{p\mid q}{1\over(p^k-1)^2} \, .
\ee
  This measure can also be represented as
\be\label{diff_sum_F}
  \xi(k)\sum_\Cond{d=1\CR d{\rm\ squarefree}}^\infty\;
  \biggl(\prod_{p\mid d}{1\over p^{2k}-2p^k}\biggr)\omega^{}_{\ZZ/d^k} \, ,
\ee
  a weak*-convergent sum of Dirac combs.
\end{theorem}
{\sc Proof\/}: As in the proof of Theorem~\ref{diff_thm}, it follows from
(\ref{ACkf}) that
\be
  \gamma\;=\;\xi(k)\sum_\Cond{d=1\CR d{\rm\ squarefree}}^\infty\;{1\over d^k}
  \prod_{p\mid d}\left(1-{2\over p^k}\right)^{-1}\omega^{}_{d^k\ZZ} \, ,
\ee
where $\gamma$ is the autocorrelation of $F\/$ and this sum of tempered
distributions is convergent in the $\|\cdot\|^{}_\loc$-topology for $k\ge2$.
Its term-by-term {}Fourier transform is
\be \label{thissum}
  \xi(k)\sum_\Cond{d=1\CR d{\rm\ squarefree}}^\infty\;{1\over d^{2k}}
  \prod_{p\mid d}\left(1-{2\over p^k}\right)^{-1}\omega^{}_{\ZZ/d^k} \, ,
\ee
where the sum of the local norms of the terms is convergent, so by
Lemma~\ref{Mtest} (\ref{thissum}) weak*-converges to a translation bounded pure
point measure equal to the pointwise sum of its terms.  Since the {}Fourier
transform operator is weak*-continuous this pure point distribution is the
diffraction spectrum of $F$.

This establishes (\ref{diff_sum_F}). To obtain (\ref{diff_amp_F}) we note
that the terms in (\ref{diff_sum_F}) that contribute to the intensity at
$a/q\/$ (where $q$ has to be $(k+1)$-free) are those with $d=m q^*$, 
where $q^*$ is the square-free kernel of $q$ and $m\in\ZZ^+$ is square-free 
and prime to $q\/$. Thus the intensity at $a/q\/$ is
\be
\xi(k)\prod_{p\mid q}{1\over p^{2k}-2p^k}
\sum_\Cond{m=1\CR m{\rm\ squarefree}\CR(m,q)=1}^\infty\;\;
\prod_{p\mid m}{1\over p^{2k}-2p^k}.
\ee
This is (\ref{point_amp}) without the factor ${\cal D}^2$ and
with $n\/$ replaced by $k$, so reduces to (\ref{diff_amp_F}). \hfill $\square$

\subsection*{Further connections and directions}

Above, we have emphasized that the sets of visible lattice points $V_\Gamma$
and the set of $k\/$-th-power-free numbers $F_k$ differ from any regular
model set (see \cite{Meyer,Moody,Martin1,Schlottmann} for definitions and properties)
by a set of positive density, suggesting that they cannot be obtained from
the cut-and-project construction in any natural way.  However, there is a
way of obtaining these sets by cut-and-project using the rational adeles
instead of Euclidean space as the hyperspace (or embedding space) and using 
windows which, although they have empty interior, are quite natural sets in this 
context. {}From this point of view, $V_\Gamma$ and $F_k$ are ``super-singular" 
model sets. This comes about because these sets are the result of sieving
over primes: $F_k$, for example, is what remains of $\ZZ$ after removing the zero 
residue class mod $p^k$ for each $p$.

The {\em cut-and-project construction} can be pictured like this.
\be \label{cutandproject}
  \begin{array}{ccccc}
   B & \stackrel{\pi}{\longleftarrow} & A &
       \stackrel{\pi^{}_\intl\;}{\longrightarrow} & C\\
   && \cup \\
   && L
  \end{array}
\ee
Usually the {\em hyperspace} $A\/$ is $\RR^{n+m}$, with the {\em physical space}
$B=\RR^n$ and the {\em internal space} $C=\RR^m$ being complementary subspaces
having $\pi$ and $\pi^{}_{\intl}$ as the associated projections with kernels $C\/$
and $B\/$ respectively, and $L$ being a lattice in $A\/$ whose images in
$B\/$ and $C\/$ are dense.  A bounded acceptance domain or {\em window} 
$\,\Omega\subset C\/$ is chosen and the {\em model set} in $B\/$ is
\be
  \Lambda \; = \; \Lambda(\Omega) \; := \; 
  \{\pi(\bm x)\mid \bm x\in L \, , \, \pi^{}_\intl(\bm x)\in \Omega \} \, .
\ee
Then $\Lambda$ is a uniformly discrete set, and also relatively dense if
$\Omega$ has non-empty interior. When $B\/$ and $C\/$ are orthogonal and
$\Omega\/$ is a simple region we have
\be\label{densM}
   \dens(\Lambda) \; = \; \vol(\Omega)\dens(L)
   \; = \; \vol(\Omega)/\vol(\cF^{}_L),
\ee
where $\cF^{}_L$ is a fundamental region of $L$. This construction admits the 
generalization where $A\/$ is allowed to be an arbitrary locally 
compact Abelian group and $L$ a discrete subgroup with $A/L$ compact
(so that $L$ plays the r\^ole of a lattice), see \cite{Moody} and
references therein.  It is well known \cite{Meyer,Moody,Martin1} that, even
in this more general situation, $\Lambda(\Omega)$ is relatively dense
if $\Omega$ has non-empty interior and possesses a uniform density
if, in addition, $\partial\Omega$ has Haar measure $0$, 
see \cite{Hof2,Martin1} for details. 

To obtain the $k\/$th-power-free numbers, we take $A$
to be the ring of rational adeles $\Adel^{}_\QQ$.  Let us explain
what this is.  (For the technical details and background to $p\/$-adic
numbers and adeles see \cite{Cassels}, which deals with the case of a
general algebraic number field in place of $\QQ$.)  Given a prime $p\in\ZZ^+$
the $p\/$-adic valuation on $\QQ$ is defined, for $q\ne0$, by $|q|_p=p^{-r}$,
where $q=p^ra/b$ with $r\in\ZZ$ and $a$, $b\/$ not divisible by $p$.  This
satisfies not only the triangle inequality but the stronger inequality
$|q_1+q_2|_p\le\max\{|q_1|_p,|q_2|_p\}$.  (Such valuations are called
{\em non-Archimedean}.)  The field $\QQ_p$ of {\em $p$-adic numbers} is the
completion of $\QQ$ with respect to $|\cdot|_p$ (analogous to $\RR$ being the
completion of $\QQ$ with respect to the ordinary absolute value) and is locally
compact in the $p\/$-adic topology.  In view of the strong triangle inequality,
the $p\/$-adic numbers $a\/$ with $|a|_p\le1$ form a ring $\ZZ_p$, called the
{\em $p$-adic integers}.  (It is the closure of $\ZZ$ in $\QQ_p$.)  {}For any
$p\/$-adic number $b\not\in\ZZ_p$ the $p\/$-adic open ball $\{a\mid|a-b|_p<1\}$
is disjoint from $\ZZ_p$.  Hence $\ZZ_p$ is closed in $\QQ_p$ and, since it is
bounded, also compact.  But $\ZZ_p$ is the disjoint union of the $p\/$ open
balls $\{a\mid|a-i|_p<1\}$ with $i=1,\ldots,p$, hence is also open in $\QQ_p$.
As a locally compact Abelian group under addition, $\QQ_p$ has a Haar measure, 
unique up to
a multiplicative constant, which can be normalized so that $\vol(\ZZ_p)=1$.

Any fixed power of a non-Archimedean valuation is also a valuation, and
topologically equivalent; but up to this equivalence the $p\/$-adic valuations
and the ordinary absolute value are the {\em only} valuations on $\QQ$.
The {\em rational adele ring} $\Adel^{}_\QQ$ is the restricted direct product
with respect to the $\ZZ_p$'s of the completions of $\QQ$ with respect to
these valuations, i.e.
\be
   \Adel^{}_\QQ \; = \;
   \Bigl\{\alpha=(\alpha_\infty,(\alpha_p))  
   \mid
   \begin{array}{c}
   \alpha_\infty\in\RR , \, \alpha_p\in\QQ_p , \mbox{ and } 
   \alpha_p\in\ZZ_p \\
   \mbox{for all except finitely many $p$.}
   \end{array} \Bigr\}
\ee
and has the restricted product topology for which the sets
$O_\infty\times\prod_pO_p$, with $O_\infty$ open in $\RR$, $O_p$ open in $\QQ_p$
and $O_p=\ZZ_p$ for all except finitely many $p$, form a base of open sets.
As a restricted product of locally compact sets with respect to compact sets,
$\Adel^{}_\QQ$ is locally compact.

Now $\QQ$ embeds in $\Adel^{}_\QQ$ diagonally (i.e. each $q\in\QQ$ can be
identified with the adele all of whose components are $q\/$), and with
this identification $\QQ$ is discrete in $\Adel^{}_\QQ$ and $\Adel^{}_\QQ/\QQ$
is compact.  A fundamental region for $\QQ$ in $\Adel^{}_\QQ$ is
$[0,1]\times\prod_p\ZZ_p$, which has volume 1 in
the normalized Haar measure on $\Adel^{}_\QQ$.  This gives rise to the
interesting and natural {\em cut-and-project scheme}
\be
  \begin{array}{ccccc}
   \RR & \stackrel{\pi}{\longleftarrow} & \Adel^{}_\QQ &
         \stackrel{\pi^{}_\intl\;}{\longrightarrow} & 
         \displaystyle\prod_p{\!}^{}_{(\ZZ_p)}\,\;\QQ_p\\[-4mm]
    && \cup \\
    && \QQ
   \end{array}
\ee
with componentwise projections and $\prod{\!}_{(\ZZ_p)}$ denoting the restricted 
product.  The image $\pi(\QQ)$ is, of course, dense
in $\RR$ and the denseness of $\pi^{}_\intl(\QQ)$ in $\prod\QQ_p$ is equivalent
to the Strong Approximation Theorem \cite[\S 15]{Cassels}.  If we choose for
the window the closed and open set $\Omega:=\prod_p\ZZ_p$ we obtain the model set
$\Lambda(\Omega)=\ZZ\in\RR$.  
Since $\vol(\Omega)=1$ in the normalized Haar measure on $\prod\QQ_p$
and $\dens(\ZZ)=1$ we notice that this satisfies (\ref{densM}).  If instead we
choose $\Omega:=\prod_p(\ZZ_p\setminus p^k\ZZ_p)$ with $k\ge2$ we obtain the ``thin"
model set $\Lambda(\Omega)=F_k$.  We say ``thin" because, being not relatively dense, 
$F_k$ is not a regular model set.  The reason that Schlottmann's result does not
apply here is that $\Omega$ has empty interior: $\Omega$ contains no basic open set
because it projects to a proper subset of $\ZZ_p$ in every non-Archimedean
component, whereas basic open sets project to $\ZZ_p$ in all except finitely
many components.  Nevertheless, $\Omega$ has a positive volume in the normalized
Haar measure on $\prod\QQ_p$ given by
\be
   \vol(\Omega) \; = \; 
   \prod_p\vol(\ZZ_p\setminus p^k\ZZ_p)
   \; = \; 
   \prod_p\Bigl( 1-{1\over p^k} \Bigr)
   \; = \; {1\over\zeta(k)} \, ,
\ee
which, in view of Proposition~\ref{prop11}, agrees with (\ref{densM}). This
is astonishing and we have no explanation for it at present, since merely 
by translating $\Omega$ we can not only reduce the density of 
$\Lambda$ to 0 but can cause $\Lambda$ to vanish altogether.
Let $i:P\to\ZZ$, where $P\/$ is the set of positive primes, be any one-one
correspondence, and for each $p\in P\/$ choose $\alpha_p\in\ZZ_p$ with
$\alpha_p\equiv i(p)\mod{p^k}$.  Then the window $\Omega+\alpha\subseteq\prod\ZZ_p$,
where $\alpha=(\alpha_p)\in\prod\QQ_p$, leads to the empty set under the
cut-and-project construction because each $i\in\ZZ$ is excluded mod $p^k$
for the corresponding $p$.

The visible points $V_\Gamma$ of a lattice $\Gamma\subset\RR^n$, for $n\ge2$,
can be obtained by a similar construction, using $(\Adel^{}_\QQ)^n$, which is
topologically isomorphic to the restricted product with respect to $\ZZ_p^n$
\be
     A \; := \; \RR^n\times\prod_p{\!}^{}_{(\ZZ_p^n)}\;\,\QQ_p^n \, .
\ee
Choose a basis $\{\bm b_1,\ldots,\bm b_n\}$ of $\Gamma$ and embed $\QQ\Gamma$
in $A\/$ by
\be
  \bm x  \; = \;
     q^{}_1\bm x^{}_1+\cdots+q^{}_n\bm x^{}_n
     \mapsto(\bm x,((q^{}_1,\ldots,q^{}_n)^{}_p)) \, ,
\ee
where $(q_1,\ldots,q_n)_p=(q_1,\ldots,q_n)$ for each $p$.  (In the terminology
of commutative algebra, this establishes that $A\/$ is isomorphic to
$\Gamma\otimes_\ZZ\Adel^{}_\QQ$, the $\Adel^{}_\QQ$-ification of the $\ZZ$-module
$\Gamma$, see \cite{Jimmy} under the index entry ``ification").  Then the
image of $\QQ\Gamma$ is a discrete subgroup of $A\/$ and $A/\QQ\Gamma$ is
compact with $\cF=\cF^{}_\Gamma\times\prod\ZZ_p^n$ as a fundamental region,
where $\cF_\Gamma$ is a fundamental region of $\Gamma$ in $\RR^n$.
The volume of this fundamental region in the normalized Haar measure
on $A\/$ is $1/\dens(\Gamma)$.  This gives the cut-and-project scheme
\be
  \begin{array}{ccccc}
   \RR^n & \stackrel{\pi}{\longleftarrow} & (\Adel^{}_\QQ)^n &
         \stackrel{\pi^{}_\intl\;}{\longrightarrow} & 
         \displaystyle\prod_p{\!}^{}_{(\ZZ_p^n)}\,\;\QQ_p^n\\[-4mm]
    && \cup \\
    && \QQ\Gamma
   \end{array}
\ee
with the images of $\QQ\Gamma$ under $\pi$ and $\pi^{}_\intl$ being dense.  If we
choose for the window the closed and open set $\Omega:=\prod_p\ZZ_p^n$ we obtain
the model set $\Lambda(\Omega)=\Gamma\subset\RR^n$ which has $\vol(\Omega)=1$ 
and satisfies (\ref{densM}).  
If we choose $\Omega:=\prod_p(\ZZ_p^n\setminus p\ZZ_p^n)$ we obtain
$\Lambda(\Omega)=V_\Gamma$.  
Again, $\Omega$ has empty interior but positive volume, given by
\be
   \vol(\Omega) \; = \; \prod_p\Bigl(1-{1\over p^n}\Bigr)
   \; = \; {1\over\zeta(n)}
\ee
so that (\ref{densM}) holds by Proposition~\ref{freq1}.  Again,
translating $\Omega$ can cause the model set to vanish.  A translation
$\bm\alpha\in\prod\QQ_p^n$ that does this can be constructed as follows: let
$\bm i:P\to\ZZ^n$ be any one-one correspondence and for each $p\/$ choose
$\bm\alpha_p\equiv\bm i(p)\mod p$; then take $\bm\alpha=(\bm\alpha_p)$.

These adelic constructions are nothing more than number theoretic sieves
which exclude certain residue classes modulo a power of each prime.  The
choice of window determines the residue classes to be retained and translating
the window changes the set of residue classes without changing their number.
Clearly, many other examples of the same type can be constructed.

If $V_\Gamma$ and $F_k$ were regular model sets, we would be able to
use standard results on model sets to derive their diffractiveness
\cite{Schlottmann}, but the constructions just described are outside
the range of current diffraction results.  To the best of our knowledge,
the diffraction of adelic model sets has not been studied, which is why
we have had to derive our diffraction results from scratch.

Obviously, one further question is how relevant the examples are. Is the set
of visible lattice points an isolated example, or more of a paradigm case of
a family of examples with pure point diffraction? 
We believe the latter is true, and there is one immediate
class of examples that springs to mind. If we view the visible points of $\ZZ^2$
as the orbit of $(1,0)$ under the group $SL(2,\ZZ)$, it is clear that this
orbit, in general, will split into several pieces when $SL(2,\ZZ)$ is replaced
by one of the standard congruence subgroups. However, our methods are robust
under further congruence constraints, and so it is clear that there is a large
class of congruence subsets of the visible points which will have interesting
diffraction spectra.

We hope to report on some examples soon.

\subsection*{Summary}

We have demonstrated that the set of visible lattice points in dimensions
$n\ge2$ and the set of $k$th-power-free integers with $k\ge2$ both possess
well-defined autocorrelations and pure point diffraction spectra. We have
shown how to replace the intuitive but incomplete derivation of the spectra
by a rigorous number theoretical argument based on explicit computation of 
the autocorrelation followed by {}Fourier inversion and repeated application
of the Poisson summation formula. 

\bigskip
\subsection*{Acknowledgements}

It is our pleasure to thank Martin N.\ Huxley for several useful hints on the
literature of $k$-free integers and Martin Schlottmann for a number of clarifying
discussions. This work was supported by the German Research Council (DFG) and
the Natural Sciences and Engineering Council of Canada (NSERC).

\clearpage
\section*{Appendix}

This paper has required several calculations of densities of discrete
point sets, for which we have used natural density.  The natural density
is in fact nowhere near as natural as one would expect from its name.
Here we briefly examine the drawbacks of the natural density and also
the uniform density, and suggest a simply stated definition of density,
intermediate between natural density and uniform density, that is
applicable to sets of the kind we encounter here.  A fuller exposition
will appear elsewhere \cite{PP}.

As an aid to describing the behaviour of densities we begin by introducing
some standard measurements of measurable sets in $\RR^n$.

\subsection*{Measurements of measurable sets}

{}For a bounded, measurable set $\cal R$ in $\RR^n$ we denote by $V({\cal R})$
the volume of $\cal R$ (that is, its $n\/$-dimensional measure).  Given
$\epsilon>0$ we denote by $\partial_\epsilon\cal R$ the set of points in
$\RR^n$ whose distance from the boundary of $\cal R$ is less than $\epsilon$.
Then $\partial_\epsilon\cal R$ is open, hence measurable, so also has a volume.
We temporarily call a set $\cal R$ {\em good\/} if
\begin{itemize}
\item[(a)]$\cal R$ is bounded and measurable, and
\item[(b)]$V(\partial_\epsilon{\cal R})/\epsilon$ is bounded
as $\epsilon\to0$.
\end{itemize}
{}For good sets we can define the functions
\begin{eqnarray}
   V&=&V({\cal R}),\\
   M&=&M({\cal R}) \; = \; \sup_{\bms x\in\cal R}|\bm x|,\\
   R&=&R({\cal R}) \; = \;
   \inf_{\bms c\in\RR^n\rule{0mm}{2.8mm}}\sup_{\bms x\in\cal R}|\bm x-\bm c|,\\
   S&=&S({\cal R}) \; = \;
   \sup_{0<\epsilon\le R}\{V(\partial_\epsilon{\cal R})/2\epsilon\}.
\end{eqnarray}

The following relations are immediate from these definitions:
\be
  R \; \le \; M,\quad V \; \le \; v_nR^n,\quad 
  RS \; \ge \; 2^{n-1}V, \quad S \; \ge \; v_nR^{n-1}/2,
\ee
where $v_n$ is the volume of the unit ball in $\RR^n$.

Each of these functions is invariant under Euclidean transformations and scales
homogeneously with dilation of $\cal R$: $M\/$ and $R\/$ proportionally to the
dilation factor $\lambda$, $V\/$ proportionally to $\lambda^n$ and $S\/$
proportionally to $\lambda^{n-1}$.  They are affected by a non-singular affine
transformation $A=T+\bm t$ (where $T\/$ is linear and \bm t is a translation)
as follows:
\begin{eqnarray}
  V(A({\cal R}))&=&|\det T|\,V({\cal R}),\label{volume}\\
  M(A({\cal R}))&\le&\|T\|\,M({\cal R})+|\bm t|,\label{bound}\\
  R(A({\cal R}))&\le&\|T\|\,R({\cal R}),\\
  S(A({\cal R}))&\le&2^n|\det T|\,\|T\|^{n-1}\|T^{-1}\|^nS({\cal R}).
\label{area}
\end{eqnarray}

Of these functions, $V$, $M\/$ and $R\/$ are straightforward, representing
the volume, the maximum distance from the origin and the circumradius of
$\cal R$, respectively.  The function $S\/$ is a substitute for the surface
area of $\cal R$ and in fact is always greater than or equal to the
$(n-1)$-dimensional measure of the boundary $\partial{\cal R}$.
The condition $\epsilon\le R$ in its definition is arbitrary,
but ensures the scaling of $S\/$ under dilations mentioned above.
(Using any fixed multiple of $R\/$ would have the same effect.)

These functions extend to arbitrary measurable sets if we define $M=R=\infty$
when $\cal R$ is unbounded (in which case $V\/$ may or may not be infinite)
and $S=\infty$ when (b) fails.  We note that $V(\partial_\epsilon{\cal R})$
is infinite for every $\epsilon$ when $\cal R$ is an unbounded set of finite
volume, so $S=\infty$ for every set of finite volume that is not good.

\subsection*{Densities of discrete sets}

Let $X\/$ be a locally finite set of points in $\RR^n$, that is, every
bounded region of $\RR^n$ contains only finitely many points of $X$.
Here we compare ways of defining a density for such a set $X$.

\subsubsection*{Uniform density}

$X\/$ has {\em uniform density\/} ${\cal D}={\cal D}(X)$ if
\be\label{uniform}
  |\{\bm x\in X\mid|\bm x-\bm c|<R\}| \; = \; {\cal D}v_nR^n+o(R^n)
\quad\mbox{as $R\to\infty$}
\ee
uniformly in \bm c.

It is immediate from this definition that\\
{\em (i) uniform density is invariant under Euclidean transformations.}\\
A less immediate consequence is that if $\cal R$ is an arbitrary measurable
set in $\RR^n$ and $f(V)$ is any function of $V\/$ (no matter how slowly
increasing) that tends to infinity with $V\/$ then
\be\label{goodcount}
  |{\cal R}\cap X| \; = \; {\cal D}V({\cal R})+o\Bigl(V({\cal R})\Bigr)+
  O\Bigl(f(V({\cal R}))S({\cal R})\Bigr)\quad\mbox{as $V({\cal R})\to\infty$}.
\ee

A consequence of (\ref{goodcount}) is\\
{\em (ii) uniform density is independent of shape;}\\
that is, if $X\/$ has uniform density $\cal D$ then (\ref{uniform})
continues to hold, with the same value of $\cal D$, when the sphere
$|\bm x-\bm c|<R$ is replaced by an expanding set of any other shape and
$v_n$ is replaced by the volume of the set in the family that has $R=1$.

Another consequence of (\ref{goodcount}) (using (\ref{volume}) and
(\ref{area}) to control the error terms) is\\
{\em (iii) uniform density varies as the reciprocal of the determinant
under non-singular affine transformations.}

All crystals and most quasicrystals possess a uniform density. However,
the set $V\/$ of the visible points of a lattice $\Gamma$ does not. In fact
Proposition~6 shows that the number of points of $V\/$ in a ball of radius
$R\/$ with centre \bm 0 is $\bigl(\dens(\Gamma)/\zeta(n)\bigr)v_nR^n+o(R^n)$,
whereas Proposition~5 shows that there are arbitrarily large balls that
contain no points of $V\/$.  Similarly, the set $F_k$ of $k\/$th-power-free
integers in $\RR$ does not possess a uniform density either.  Consequently
uniform density is of no use for the questions considered in this paper---in
fact none of the sets whose densities we require possesses a uniform density.

A definition of density known as {\em van Hove density\/} is implicit in
Ch.~2 of \cite{Ruelle}, if the identical formul\ae\ (3.3) and (3.8) of
\cite[Ch.~2]{Ruelle} are regarded as defining the density of a general
potential operator $\Phi$.  This is equivalent to uniform density.
The fact that uniform density implies van Hove density follows from
(\ref{goodcount}) and the converse implication results from choosing
the van~Hove sets to be balls.

\subsubsection*{Natural density}

As defined in (10) for uniformly discrete sets, $X\/$ has
{\em natural density\/} ${\cal D}={\cal D}(X)$ if
\be
  |\{\bm x\in X\mid|\bm x|<R\}| \; = \; {\cal D}v_nR^n+o(R^n)
  \quad\mbox{as $R\to\infty$}.
\ee
This is often called {\em asymptotic density\/} in the context of
number theory.

Like uniform density, natural density is\\
{\em (i) invariant under Euclidean transformations.}\\
(Tranlation invariance is easy to derive, though not quite as immediate
as for uniform density.)  However, natural density is not independent
of shape and does not transform naturally under linear transformations.
{}For example, let $X=\{(x,y)\mid x,y\in\ZZ, |y|<2|x|\}$ be the set of
points in $\ZZ^2$ that lie in the double wedge with angle $2\tan^{-1}2$.
The natural density of $X\/$ is $(2/\pi)\tan^{-1}2=0.7048\dots\,$, but the
proportion of the integer points in the square $\{(x,y)\mid|x|,|y|<R\}$
that belong to $X\/$ is 3/4 and the proportion of the integer points in
the square $\{(x,y)\mid|x|+|y|<R\}$ that belong to $X\/$ is 2/3.  So for
this set, counting points according to the $L^2$, $L^\infty$ or $L^1$ norm
in $\RR^2$ gives three distinct answers.  The linear map $T(x,y)=(2x,y)$
transforms $X\/$ into the set of integer points with $x\/$-coordinate even
in the double wedge $\{(x,y)\mid x,y\in\ZZ, |y|<|x|\}$ with angle $\pi/2$.
This set has density $1/4$, which is not a half of $(2/\pi)\tan^{-1}2$,
showing that natural density does not transform naturally under linear maps.

{}For the sets we deal with in this paper natural density does not exhibit
such pathological behaviour.  {}For example, Proposition~6, Theorem~1 and
Theorem~2, which establish the densities of certain subsets of a lattice
$\Gamma$ in $\RR^n$, all give a factor $\dens(\Gamma)$ as the only dependency
of the density on the lattice, confirming that the density varies as the
reciprocal of the determinant under linear transformations.  However, rather
than having to verify this in each individual case, as in effect we have
done here, it would be preferable to have a definition of density with
strength intermediate between uniform density and natural density for
which non-pathological behaviour is guaranteed.

\subsubsection*{Tied density}

It is possible to define an alternative form of density, intermediate
 between uniform density and natural density, that overcomes these
 difficulties.  This depends on the idea of making uniform density less free of the
 origin by replacing $R\/$ by $M\/$ in the error term of (\ref{uniform}).
 We say that $X$ has {\em tied density} ${\cal D}={\cal D}(X)$ if
 for all $\bm{c}\in\RR^n$
\be\label{intermediate}
  |\{\bm x\in X\mid|\bm x-\bm c|<R\}| \; = \; {\cal D}v_nR^n+o(M^n)
  \quad\mbox{as $R\to\infty$},
\ee
where $M=R+|\bm c|$.

Again it is easy to see that\\
{\em (i) tied density is invariant under Euclidean transformations.}\\
In a similar way to (\ref{goodcount}) it can also be shown that if $X\/$ has
tied density $\cal D$ and $\cal R$ is an arbitrary measurable set in
$\RR^n$ then 
\be
  |{\cal R}\cap X| \; = \; {\cal D}V({\cal R})+o\Bigl(M({\cal R})^n\Bigr)+
  o\Bigl(M({\cal R})S({\cal R})\Bigr)\quad\mbox{as $V({\cal R})\to\infty$}.
\ee
This is sufficient to establish that\\
{\em (ii) tied density is independent of shape,}\\
and, with the help of (\ref{volume}) and (\ref{bound}), that\\
{\em (iii) tied density varies as the reciprocal of the determinant
under non-singu\-lar affine transformations.}

The error term $o(M^n)$ in the definition of tied density has the
effect that, while balls centred away from the origin are not excluded
from consideration, those distant from the origin by more than a few radii
have little influence.  It can be regarded as a reminder to experimentalists
that, as the radius of a probe is increased, the next position of the probe
should not be more than a few radii from the previous position. There is,
nevertheless, a great deal of arbitrariness in the definition of tied
density.  The core idea is to make the error term in (\ref{intermediate})
depend on $\bm c$ as well as on $R$, but any kind of dependency on $\bm c$
(provided that for every $R\/$ the error term tends to infinity with $\bm c$)
would have a similar effect.  The faster the rate of increase of the error
term with $\bm c$ the more widely applicable the density definition becomes.
In an effort to conceal this arbitrariness we have chosen a very simple
dependency on $\bm c$.

It is clear from these definitions that the existence of uniform density
for a set $X\/$ implies the existence of tied density for $X$,
which in turn implies the existence of natural density for $X$, and that
each of these densities has a unique value when it is defined.

To establish the existence of tied density for the sets studied in
this paper we would need to count points $\bm x$ with $|\bm x-\bm c|<R\/$
instead of $|\bm x|<R$.  As a result, though the $\bm x$ and $\bm y$ variables
in sums still have range $R\/$ (but no longer centred at the origin), the $l\/$
and $m\/$ variables have range $R+|\bm c|=M$.  The effect is that main terms
remain the same but some $R\/$'s are replaced by $M\/$'s in error terms.
The new error terms for the results of this paper when the ball of radius
$R\/$ is centred away from the origin are as follows:
\begin{center}
\begin{tabular}{ll}
Proposition~6 ($n\ge3$)&$O(1/R)+O(M/R^n)$,\\
Proposition~6 ($n=2$)&$O(\log M/R)+O(M/R^2)$,\\
Theorem~1 ($n\ge3$)&$O(1/R)+O(M^2/R^n)$,\\
Theorem~2 ($n=2$)&$O(\sqrt M/R)+O(M/R^2)$,\\
Proposition~11 ($n=1$)&$O(M^{1/k}/R)$,\\
Theorem~4 ($n=1$)&$O(R^{1/k}M^{1/k}/RM^{1/k^2})+O(M^{1/k}/R)$.
\end{tabular}
\end{center}
Each of these error terms, when multiplied by $R^n$, is well within the
error estimate $o(M^n)$ required by the definition (\ref{intermediate}) of
tied density.  Consequently we could have worked with tied
densities instead of natural densities throughout.  One concrete advantage
of such an approach would have been that we could then have worked with
the visible points of cubic lattices only, then used property (iii) of
tied density to transfer the results to other lattices.


\clearpage
\begin{thebibliography}{99}
\small

\bibitem{Apostol}
T.~M.~Apostol, 
{\em Introduction to Analytic Number Theory},
Springer, New York (1976).

\bibitem{Arga}
L.~Argabright and J.~Gil de Lamadrid,
{\em {}Fourier analysis of unbounded measures on locally
compact Abelian groups}, Memoirs of the American Mathematical
Society, vol.\ 145, AMS, Providence (1974).

\bibitem{BGW}
M.~Baake, U.~Grimm and D.~H.~Warrington,
``Some remarks on the visible points of a lattice'',
{\em J.\ Phys.} {\bf A27} (1994) 2669--74 and
5041 (Erratum); see also math-ph/9903046.

\bibitem{Berberian}
S.~K.~Berberian, 
{\em Measure and Integration},
Macmillan, New York (1965), reprint: Chelsea, New York (1970).

\bibitem{Berend}
D.~Berend, in preparation.

\bibitem{Cassels}
J.~W.~S.~Cassels, 
``Global {}Fields'', in: {\em Algebraic Number Theory},
eds.\ J.~W.~S.~Cassels and  A.~Fr\"ohlich, Academic Press, London (1967),
pp.\ 42--84.

\bibitem{Cowley}
J.~M.~Cowley, 
{\em Diffraction Physics}, 3rd ed., North Holland, Amsterdam (1995).

\bibitem{Dieu}
J.~Dieudonn\'{e},
{\em Treatise on Analysis}, vol.\ II, Academic Press, New York (1970).

\bibitem{Grimmett}
G.~Grimmett, ``Large deviations in the random sieve'',
{\em Math.\ Proc.\ Cambridge Phil.\ Soc.} {\bf 121} (1997) 519--30.

\bibitem{GL}
P.~M.~Gruber and C.~G.~Lekkerkerker, 
{\em Geometry of Numbers}, 
2nd ed., North Holland, Amsterdam (1987).

\bibitem{HR}
E.~Hewitt and K.~A.~Ross,
{\em Abstract Harmonic Analysis I}, 2nd ed.,
Springer, New York (1979).

\bibitem{Hof}
A.~Hof,
``On diffraction by aperiodic structures'',
{\em Commun.\ Math. Phys.} {\bf 169} (1995) 25--43.

\bibitem{Hof2}
A.~Hof,
``Diffraction by aperiodic structures'', in:
{\em The Mathematics of Long-Range Aperiodic Order},
ed.\ R.~V.~Moody, NATO ASI Series C 489, Kluwer, Dordrecht (1997), 
pp.\ 239--68.

\bibitem{Hua}
L.~K.~Hua, 
{\em Introduction to Number Theory}, Springer, Berlin (1982).

\bibitem{Huxley}
M.~N.~Huxley, 
``Moments of differences between square-free numbers'', in:
{\em Sieve Methods, Exponential Sums, and their Applications 
in Number Theory}, eds.\ G.\ R.\ H.\ Greaves, G.\ Harman and M.\ N.\ Huxley,
LMS 237, Cambridge University Press, Cambridge (1997), pp.\ 187--204.

\bibitem{Jimmy}
J.~T.~Knight,
{\em Commutative Algebra},
Cambridge University Press, Cambridge (1971).

\bibitem{Lagarias}
J.~C.~Lagarias, 
``Meyer's concept of quasicrystal and quasiregular sets'',
{\em Commun.\ Math.\ Phys.} {\bf 179} (1996) 365--76.

\bibitem{Meyer}
Y.~Meyer,
{\em Algebraic Numbers and Harmonic Analysis},
North-Holland, Amsterdam (1972).

\bibitem{Mirsky}
L.~Mirsky,
``Arithmetical pattern problem relating to divisibility by $r$th powers'',
{\em Proc.\ London Math.\ Soc.} {\bf 50} (1949) 497--508.

\bibitem{Moody}
R.~V.~Moody,
``Meyer sets and their duals'', in:
{\em The Mathematics of Long-Range Aperiodic Order},
ed.\ R.~V.~Moody, NATO ASI Series C 489, Kluwer, Dordrecht (1997), 
pp.\ 403--41.

\bibitem{Mosseri}
R.~Mosseri, 
``Visible points in a lattice'',
{\em J.\ Phys.} {\bf A25} (1992) L25--9.


\bibitem{Pillai}
S.~S.~Pillai,
``On sets of square-free numbers'', 
{\em J.\ Indian Math.\ Soc.\ (New Series)} {\bf 2} (1936) 116--8.

\bibitem{PP}
P.~A.~B.~Pleasants, in preparation.

\bibitem{RS}
M.~Reed and B.~Simon,
{\em Methods of Modern Mathematical Physics. I: {}Functional Analysis},
2nd ed., Academic Press, San Diego (1980).

\bibitem{Rudin}
W.~Rudin,
{\em {}Functional Analysis}, 2nd ed.,
McGraw-Hill, New York (1991).

\bibitem{Ruelle}
D.~Ruelle,
{\em Statistical Mechanics: Rigorous Results},
Benjamin, Reading (1969); reprint: Addison-Wesley, Redwood City (1989).

\bibitem{Martin1}
M.~Schlottmann,
``Cut-and-project sets in locally compact Abelian groups'', in:
{\em Quasicrystals and Discrete Geometry}, ed.\ J.\ Patera,
Fields Institute Monographs, vol.\ 10, AMS, Providence (1998), pp.\ 247--64.

\bibitem{Schlottmann}
M.~Schlottmann, 
``Generalized model sets and dynamical systems'',
to appear in: {\em Directions in Mathematical Quasicrystals},
eds.\ M.~Baake and R.~V.~Moody,
CRM Monograph Series, AMS, Providence (2000), in press.

\bibitem{Schroeder}
M.~R.~Schroeder, 
``A simple function and its {}Fourier transform'',
{\em Mathem.\ Intelligencer} {\bf 4} (1982) 158--61,
and: {\em Number Theory in Science and Communication},
3rd ed., Springer, Berlin (1997).

\bibitem{Schwartz}
L.~Schwartz, 
{\em Th\'eorie des Distributions}, 3rd ed., Hermann, Paris (1998).

\end{thebibliography}
\end{document}